%% file: main.tex
\documentclass{amsart}

\usepackage[foot]{amsaddr}
\usepackage{quiver}
\usepackage{amsthm}
\usepackage{amsmath}
\usepackage{mathrsfs}
\usepackage{graphicx} 
\usepackage{faktor} 
\usepackage{MnSymbol}       
\usepackage{stmaryrd}
 \usepackage{relsize}
 \usepackage{hyperref}
 \usepackage{musicography}
\usepackage{euscript,enumitem}
\usepackage{bbm}
\usepackage{geometry}
\usepackage[all]{xy}
\usepackage[utf8]{inputenc}
\usepackage[english]{babel}
\usepackage{amsrefs}

\usepackage{standalone} 
\usepackage{tikzscale} 
\usetikzlibrary{backgrounds,fit,positioning}
\usetikzlibrary{shapes.geometric}
\usetikzlibrary{positioning}

\input settings.tex

\title{Boardman-Vogt tensor product and Wreath product of operadic categories}
\author{Daria Pavlova}
\email{daria.pavlova524@student.cuni.cz}
\address{Faculty of Mathematics and Physics, Charles University \& Institute of Mathematics, Czech Academy of Sciences}


\begin{document}

\begin{abstract}
We introduce the wreath product of operadic categories and show that it produces a monoidal structure on $\CatOp$, the category of strict operadic categories. We relate the wreath product to the well known Boardman-Vogt tensor product of operads by interpreting the category of symmetric set-valued colored unital operads as a reflexive subcategory of $\CatOp$. Our main
result is that the reflection is a strong monoidal functor. In particular, the Boardman–Vogt tensor product of two colored operads in $\Set$ is isomorphic to an operad that is induced by the wreath product of suitable operadic categories.
\end{abstract}
 
\maketitle
\setcounter{tocdepth}{1}
\tableofcontents

\setlength{\baselineskip}{15pt plus 1pt}
\section*{Introduction}
\input{intro}

\paragraph{\textbf{Acknowledgment}:} The author thanks Michael Batanin for suggesting this project and both Batanin and Martin Markl, the author’s advisor, for their guidance.
 
I would also like to thank my friends and colleagues, Maroš Grego and Dominik Trnka, for their valuable contributions to the discussions on this topic.

During this project, the author was supported by  Praemium Academiæ of M. Markl, by Charles University Research Center program
No. UNCE/24/SCI/022, the project SVV-2025-260837, by the GA UK project No.~433125 and by RVO: 67985840.

\section{Symmetric operads as a reflective subcategory} \label{AI}
\input{op_vs_opcat}

\section{Wreath product and the Boardman-Vogt product} \label{wreath}
\input{wreath_prod}

\section{Relation to other wreath products}
\input{relation}

\setlength{\baselineskip}{15pt}

\bibliography{bibliography}

\end{document}

%% file: settings.tex
\theoremstyle{definition}
\newtheorem{thm}{Theorem}

\newtheorem{prop}[thm]{Proposition}
\newtheorem{cor}[thm]{Corollary}
\newtheorem{obs}[thm]{Observation}
\newtheorem{notat}[thm]{Notation}

\theoremstyle{definition}
\newtheorem{defn}[thm]{Definition}
\newtheorem{example}[thm]{Example}
\newtheorem{rem}[thm]{Remark}

\newcommand{\opcat}[1]{\texttt{#1}}
\newcommand{\oper}[1]{\EuScript #1}

\newcommand{\Fin}{\texttt{Fin}}
\newcommand{\Dalg}{\Delta_{\mathrm{alg}}}
\def\Fib{{\textit{Fib}}}
\newcommand{\Set}{\textit{Set}}

\newcommand{\CatOp}{\opcat{OpCat}}
\newcommand{\SOp}{\texttt{SOp}}
\newcommand{\Op}{\texttt{Op}}
\newcommand{\Opc}[2]{\texttt{Op}^{#1}_{#2}}
\newcommand{\V}{\mathscr{V}}
\renewcommand{\O}{\opcat{O}}
\newcommand{\wprod}[2]{{\displaystyle #1 ~\mathlarger{\wr}~ #2}}

\newcommand{\Bq}{\texttt{Bq}}
\newcommand{\BV}{\otimes_{BV}}

\newcommand{\I}{\mathbb{I}}
\newcommand{\A}{\mathbb{A}}

\newcommand{\OwR}{\opcat{O} \wr \opcat{R}}

\newcommand{\card}[1]{| #1 |}
\newcommand{\termop}{\zeta}

\def\ccdots{\mathinner{\cdotp\mkern-4mu\cdotp\mkern-4mu\cdotp}}
\def\frC{{\mathfrak C}}
\def\frD{{\mathfrak D}}
\renewcommand\to{\longrightarrow}
\renewcommand\mapsto{\longmapsto}

\newcommand{\AR}{\texttt{Arity}}
\newcommand{\Arleft}{\texttt{A}}

\def\mytreesq[#1,#2,#3,#4,#5,#6,#7,#8]{\tikz[baseline=.1ex]{
\tikzset{square/.style={regular polygon,regular polygon sides=4}}
\tikzset{bl/.style = {shape=circle,draw=black,inner sep=1.5pt,font=\tiny}}
\tikzset{vert/.style = {shape=circle,fill=black,inner sep=1.5pt,font=\tiny}}
    \node at (0,-1.5) [square,draw] (root) {$S$};
    \node [#3] (row1) [label={[label distance=.2cm]0:#1}] at (0,-.5) {};
    \node (dots1) at (0,.5) {$\cdots$};
    \node [#4] (row2_1) [label=left:#2] at (-2, .5) {};
    \node [#4] (row2_2) [label=right:#2] at (2, .5) {};
    \node at (-2.8,1.5) [draw,black,rectangle] (row3_1) {#5};
    \node (dots1) at (-2,1.5) {$\cdots$};
    \node at (-1.2,1.5) [draw,black,rectangle] (row3_2) {#6};
    \node at (1.2,1.5) [draw,black,rectangle] (row3_3) {#7};
    \node (dots1) at (2,1.5) {$\cdots$};
    \node at (2.8,1.5) [draw,black,rectangle] (row3_4) {#8};

    \draw (root) -- (row1);
    \draw (row1) -- (row2_1);
    \draw (row1) -- (row2_2);
    \draw (row2_1) -- (row3_1);
    \draw (row2_1) -- (row3_2);
    \draw (row2_2) -- (row3_3);
    \draw (row2_2) -- (row3_4);
}}

%% file: intro.tex
M.~Boardman and R.~Vogt defined an associative and commutative product 
\hbox{$\oper{P} \BV \oper{Q}$} of operads $\oper{P}$ and $\oper{Q}$ ~\cite{boardman-vogt}. This product
is characterized by the property that \hbox{$(\oper{P} \BV \oper{Q})$}-algebras are the same as \hbox{$\oper{P}$-algebras} in the category of $\oper{Q}$-algebras, or equivalently, $\oper{Q}$-algebras in the category of $\oper{P}$-algebras. One way to give the Boardman-Vogt product a constructive definition is to say that \hbox{$\oper{P} \BV \oper{Q}$} is the coproduct of operads \hbox{$\oper{P} \coprod \oper{Q}$} modulo the \textit{interchange} relation \cite{Brinkmeier2010}:

\begin{center}
\includegraphics[width=0.7\linewidth]{interchange_intro.tex}
\end{center}


The interchange, however, creates an internal structure that is difficult to handle explicitly. It is a known fact that the Boardman-Vogt tensor product does not preserve homotopy equivalences. Moreover, the seemingly intuitive statement that the Boardman-Vogt tensor product of $E_n$-operads is additive, namely that $E_n \otimes E_m$ has the homotopy type of $E_{n+m}$, does not hold in general. One of the most well-known counterexamples is the tensor product \hbox{$(\oper{Ass} \BV \oper{Ass}) \cong Comm$}, where $Ass$ is the associative operad (an $E_1$-operad) and $Comm$ is the commutative operad (an $E_{\infty}$-operad). A version of the additivity property was proven by Dunn \cite{dunn}, concretely that an $n$-fold tensor product of \textit{the little 1-disk operad} is an $E_n$-operad in 1988. Then, in 2000, M.~Brinkmeier proved that the tensor product of little $n-$disks and $m-$disks operads has the homotopy type of the little $(m+n)$-disks operad \cite{Brinkmeier2010} (this proof was later revised and shortened by M. Barata and I. Moerdijk in 2024 \cite{baratamoerdijk}). In 2007, M. Brun, Z. Fiedorowicz, and R. Vogt published a proof that the tensor product of the operad $Ass$ and the little $n$-cubes operad is an $E_{n+1}$-operad \cite{brun}. And finally, Z. Fiedorowicz and R. Vogt proved in 2015 that given a cofibrant $E_n$-operad and a cofibrant $E_m$-operad, their Boardman-Vogt tensor product $E_n \otimes E_m$ is indeed an  $E_{n+m}$-operad~\cite{FIEDOROWICZ2015421}.

The purpose of the present text is to incorporate the Boardman-Vogt product into the paradigm of operadic categories, defined by M. Batanin and M. Markl in \cite{opcat}, which will hopefully provide a better understanding of the additivity theorem in our future work. We introduce a monoidal structure on the category of operadic categories $\CatOp$ with the monoidal product called \textit{the wreath product} of operadic categories. We then show that the category $\SOp$ of colored symmetric unital operads valued in sets is a reflexive subcategory of $\CatOp$. Moreover, the left adjoint $\A:~\CatOp~\to~\SOp$ to the inclusion $\I:~\SOp~\hookrightarrow~\CatOp$ is a strong monoidal functor, where $\SOp$ is endowed with the monoidal structure given by the Boardman-Vogt tensor product. This property of the functor $\A$ implies the main result of this paper, i.e. that given $\oper{P}, \oper{Q} \in \SOp$ it holds
\begin{equation}\label{eq-main}
    \oper{P} \BV \oper{Q} \cong \A(\I\oper{P} \wr \I\oper{Q}).
\end{equation}

So far, this result only applies to operads valued in the category $\Set$ of sets and arbitrary maps between sets, since the inclusion $\I: \SOp \hookrightarrow \CatOp$ is given by the operadic Grothendieck construction, which is originally defined only for non-enriched operads. However, we believe that the result we obtain is still important. When considering topological operads and their homotopy equivalences, for example, the underlying set structure may still be of interest. This can be seen in the proofs by Brinkmeier, Barata, and Moerdijk mentioned above. There are relevant generalizations of the Grothendieck construction for enriched cases that we hope to explore further in future work, namely one for non-symmetric operads with values in categories by D. Trnka \cite{trnka} and one for categories enriched over monoidal $\V$ with some further requirements by J.~Beardsley and L.~Ze~Wong \cite{enriched-gr}.

Similar kinds of products have been defined for categories with some additional structure. C.~Berger introduced a  product called the \textit{categorical} wreath product in \cite{BERGER2007230} for categories equipped with a functor to the category $\Gamma$, which is the opposite of the category of finite pointed sets. The application of the latter to the category $\Delta$ of finite non-empty ordinals and order-preserving maps produces the category $\Theta_n$, which then characterizes $n$-fold loop spaces. Berger, moreover, observes that the canonical \textit{assembly} functor corresponding to the iterated categorical wreath product on $\Delta$ (i.e., the category $\Theta_n$) is dual (up to a certain restriction) to the \textit{cardinality} functor corresponding to the iterated operadic wreath product on $\Dalg$ (i.e., the category $\Omega_n$). We expand on this observation and show a kind of duality between the categorical and the operadic wreath products. 

Another kind of wreath product was introduced by C.~Barwick for \textit{operator categories}, introduced in \cite{barwick-paper}. Operator categories can be viewed as a particular type of operadic category, and the wreath product of the latter is a natural generalization of the wreath product of the former. Moreover, the wreath product of \textit{operator} categories induces a generalized version of the Boardman-Vogt tensor product of operads over \textit{operator} categories. 

In future work, we wish to define a suitable product on the category of all operads over various operadic categories induced by the wreath product of operadic categories and study its relation to the Boardman-Vogt tensor product. We are motivated by a noteworthy property of the wreath product of operadic categories, which was also observed by Berger and Barwick in their respective wreath products: when applied iteratively to the category of finite (including empty) ordinals and order-preserving maps, the construction yields the category $\Omega_n$ of Batanin's $n$-trees \cite[Sec.~3, Ex.~8]{BATANIN199839}. Our aim in the future is to study the induced product of $n$-operads and its relation to the tensor product of $E_n$-operads. This paper lays some groundwork for this.

\paragraph{\textbf{Organization of the paper.}}
In Section \ref{AI}, we recall the relevant definitions of the theory of operadic categories and make a few useful observations regarding the canonical Arity functor. Then we establish the framework through which the main result is shown. This involves demonstrating an adjunction $\xymatrix{\A:\CatOp\ar[r]<1ex> & \SOp:\I\ar[l]<1ex>_(0.44){\upvdash}}$ between the category $\CatOp$ of strict operadic categories and strict operadic functors, and the category $\SOp$ of colored symmetric operads in $\Set$ and morphisms of operads. We show that $\SOp$ is a reflective subcategory of $\CatOp$ in Theorem \ref{thm-adj-AI}, which is a result of independent interest.

In Section 2, we introduce a monoidal structure on the category $\CatOp$ given by the \textit{wreath} product of operadic categories. We then show that the aforementioned functor $\A: \CatOp \to \SOp$ is strong monoidal in Theorem \ref{thm-mon-fun}, which then implies the isomorphism (\ref{eq-main}).

In Section 3, we recall C.~Berger's categorical wreath product and C.~Barwick's wreath product of operator categories and overview their correspondence to the wreath product of operadic categories introduced in this paper.

\paragraph{\textbf{Conventions}} Given an operad $\oper{P}$ in the sense of May \cite{may}, we always denote the composition maps by $\gamma_{\oper{P}}$. Given an operad $\oper{Q}$ in the sense of Batanin and Markl \cite{opcat}, we always denote the composition maps by $\mu_{\oper{Q}}$. We omit the subscript when the operad is clear from the context.


%% file: op_vs_opcat.tex
We recall some basic definitions from \cite{opcat} and make a few useful observations about the canonical Arity functor. We use these observations to show that colored operads form a reflective subcategory of operadic categories. For brevity, we use the terms \textit{operadic category} and \textit{operadic functor} to refer to what have been defined as a \textit{strict operadic category} and a \textit{strict operadic functor} in \cite{opcat}.

Let $\Fin$ be the skeletal category of finite sets. The objects of
this category are linearly ordered sets $\underline{n} =\{1 \leq \ldots \leq n\},
n\in \mathbb{N}.$ We sometimes omit the underline notation and simply write $n$ for the respective linearly-ordered set. Morphisms are arbitrary (not necessarily order-preserving) maps between the underlying sets. We define the $i$-th fiber $f^{-1}(i)$ of a morphism $f: T \to S$, $i\in S$, as the pullback of $f$ along the map $\ulcorner{i}\urcorner: \underline{1} \to S$ which picks up the element $i$.

Throughout the paper we identify elements of the fiber given as a pullback in $\Fin$ with their  order-preserving inclusion to the preimage. 

Any commutative triangle
\[\begin{tikzcd}[ampersand replacement=\&,cramped, column sep=tiny]
	T \&\& S \\
	\& R
	\arrow["f", from=1-1, to=1-3]
	\arrow["h"', from=1-1, to=2-2]
	\arrow["g", from=1-3, to=2-2]
\end{tikzcd}\]
in $\Fin$ induces a map $f_i: h^{-1}(i) \to g^{-1}(i)$, for each $i \in R$. Moreover, this assignment is functorial, and the equality $f^{-1}(j) = f^{-1}_{g(j)}(j)$ holds for any $j \in S$. The above structure on the category $\Fin$ motivates the structure required for an operadic category.

\begin{defn}
An \textit{operadic category} $\O$ is a category equipped with a \textit{cardinality} functor \[\card{-}: \O \to \Fin\] that has the following properties. We require that each connected component of $\O$ has a \textit{chosen local terminal object} $U_c$, $c \in \pi_0(\O)$. We also assume that for every $f: T \to S$ in $\O$ and every element $i \in \card{S}$, there is an object $f^{-1}(i)$, which we will call the \textit{$i$-th fiber of $f$}, such that $\card{f^{-1}(i)} = \card{f}^{-1}(i)$. We use the notation $f^{-1}(i) \rhd T \xrightarrow{~f~}S$ to indicate the fibers. We also require that 

\begin{itemize}  
\item[(AX1)] \label{ax1} For any $c\in \pi_0(\opcat{O})$, $|U_c| = 1$.
\end{itemize}

A trivial morphism $f:T\to S$ in $\opcat{O}$ is a morphism such that, for each $i\in
 |S|$,  
$f^{-1}(i) = U_{d_i}$ for some $d_i\in \pi_0(\opcat{O}).$

The remaining axioms for a strict operadic category are:
\begin{itemize}
\item[(AX2)]  \label{ax2} The identity morphism $id:T\to T$ is trivial for any $T\in \opcat{O};$

\item[(AX3)] For any commutative diagram in $\opcat{O}$
  \begin{equation*}\label{diag-axiom3}
\begin{tikzcd}[ampersand replacement=\&, column sep=tiny]
	T \&\& S \\
	\& R
	\arrow["f", from=1-1, to=1-3]
	\arrow["h"', from=1-1, to=2-2]
	\arrow["g", from=1-3, to=2-2]
\end{tikzcd}
\end{equation*}
  and every $i\in |R|$, one is given a map
  \[
  f_i: h^{-1}(i)\to g^{-1}(i)
  \]
  such that $|f_i|: |h^{-1}(i)|\to |g^{-1}(i)|$ is the map
  $|h|^{-1}(i)\to |g|^{-1}(i)$ of sets induced by
  \begin{equation*}
\begin{tikzcd}[ampersand replacement=\&, column sep=tiny]
	{\card{T}} \&\& {\card{S}} \\
	\& {\card{R}}
	\arrow["{\card{f}}", from=1-1, to=1-3]
	\arrow["{\card{h}}"', from=1-1, to=2-2]
	\arrow["{\card{g}}", from=1-3, to=2-2]
\end{tikzcd}
  \end{equation*}
  We moreover require that this assignment forms a functor ${\Fib}_i: \opcat{O}/R \to \opcat{O}$. 
  
  \item[(AX3a)]If
  $R=U_c$, the functor ${\Fib}_1$ is required to be the domain functor
  $\opcat{O}/R \to \opcat{O}.$
  
 \item[(AX4)]
In the situation of (AX3), for any $j \in |S|$ and $i = \card{g}(j)$, there is an equality of fiber objects
\begin{equation*}
\label{axiom-4}
\begin{tikzcd}[ampersand replacement=\&,cramped, column sep=tiny]
	{f^{-1}_i(j)} \& {h^{-1}(i)} \&\& {g^{-1}(i)} \\
	{f^{-1}(j)} \& T \&\& S \\
	\&\& R
	\arrow["\rhd"{marking, allow upside down, pos=0.3}, draw=none, from=1-1, to=1-2]
	\arrow["{=}"{marking}, draw=none, from=1-1, to=2-1]
	\arrow["{f_i}", from=1-2, to=1-4]
	\arrow["\rhd"{marking, allow upside down}, draw=none, from=1-2, to=2-2]
	\arrow["\rhd"{marking, allow upside down}, draw=none, from=1-4, to=2-4]
	\arrow["\rhd"{marking, allow upside down, pos=0.2}, draw=none, from=2-1, to=2-2]
	\arrow["f", from=2-2, to=2-4]
	\arrow["h"', from=2-2, to=3-3]
	\arrow["g", from=2-4, to=3-3]
\end{tikzcd}
\end{equation*}

 \item[(AX5)] Let 
\begin{equation*}
\begin{tikzcd}[ampersand replacement=\&]
	\& S \\
	T \& { } \& Q \\
	\& R
	\arrow["a", from=1-2, to=2-3]
	\arrow["g"{pos=0.2}, from=1-2, to=3-2]
	\arrow["f", from=2-1, to=1-2]
	\arrow["b"{pos=0.8}, no head, from=2-1, to=2-2]
	\arrow["h"', from=2-1, to=3-2]
	\arrow[from=2-2, to=2-3]
	\arrow["c", from=2-3, to=3-2]
\end{tikzcd}
\end{equation*}
be a commutative diagram in $\opcat{O}$ and let $i \in |Q|, j = |c|(i)$.
Then by (AX3) the diagram
\begin{equation*}
\begin{tikzcd}[ampersand replacement=\&, column sep=tiny]
	{h^{-1}(i)} \&\& {g^{-1}(i)} \\
	\& {c^{-1}(i)}
	\arrow["{f_i}", from=1-1, to=1-3]
	\arrow["{b_i}"', from=1-1, to=2-2]
	\arrow["{a_i}", from=1-3, to=2-2]
\end{tikzcd}
\end{equation*}
commutes, so it induces a morphism $(f_i)_{j}: b_i^{-1}(j)\to a_i^{-1}(j).$ By axiom (AX4) we have $$a^{-1}(j)=a_k^{-1}(j) \ \mbox{and} \ b^{-1}(j)=b_i^{-1}(j).$$
We  then require the equality \begin{equation*}
    f_i = (f_i)_{j}.
\end{equation*}
\end{itemize}
We will also assume that the set $\pi_0(\opcat{O})$ of connected components
is {\em small\/} with respect to a sufficiently big ambient universe.
\end{defn}

An \textit{operadic functor} between two operadic categories is a functor $F: \O \to \opcat{R}$ that commutes with the cardinality functor,  preserves fibers, local terminal objects, induced morphisms, and equalities required by the axioms of operadic categories. This defines the category $\CatOp$ of operadic categories and operadic functors. 

\begin{example} \label{ex-delta}
    The category $\Fin$, equipped with the identity functor as its cardinality and with fibers defined via pullbacks, is evidently an operadic category. Its subcategory $\Dalg$, which has the same objects $\underline{n} = \{1, \ldots, n\}, n \geq 0$, including the object $\underline{0} = \emptyset$, but only order-preserving maps, is likewise an operadic category when endowed with the inclusion $\Dalg \hookrightarrow \Fin$ as its cardinality functor and with fibers given by pullbacks.

    In contrast, the classical simplex category $\Delta$, whose objects are denoted by $[n] = \{0 \leq \cdots \leq n\}$ for $n \geq 0$, does not contain an empty ordinal. Consequently, it cannot be endowed with the structure of an operadic category via the inclusion $[n] \mapsto \underline{n+1}$ as its cardinality functor, since fiber objects for non-surjective morphisms would not exist.

    We emphasize the difference in the notation of the objects of $\Dalg$ and $\Delta$.
\end{example}

\begin{example}
Let $\frC$ be a set. A $\frC$-bouquet is a map $b: \underline{k} +1 \to \frC$, where $\underline{k} \in \Fin$. In other words, a $\frC$-bouquet is an ordered $(k+1)$-tuple $(c_1, \ldots, c_k; c)$ of elements of $\frC$. It can be viewed as a planar corolla in which all edges, including the root, are colored by elements of $\frC$. 

\[
\begin{tikzpicture}
    \draw (0,0) -- (0,0.7); 
    \draw (0,0.7) -- (-0.3,1.3); 
    \draw (0,0.7) -- (0.3,1.3); 
    \draw (0,0.7) -- (-0.8,1.3); 
    \draw (0,0.7) -- (0.8,1.3); 
    
    \node (c) at (0, -0.2) {$c$};
    \node (c1) at (-0.8, 1.5) {$c_1$};
    \node (c2) at (-0.3, 1.5) {$c_2$};
    \node (c3) at (0.3, 1.5) {$c_3$};
    \node (c4) at (0.8, 1.5) {$c_4$};
\end{tikzpicture}
\]

The extra color $c \in \frC$ is called the root color. The finite set $\underline{k}$ is the underlying set of the bouquet~$b$. A map of $\frC$-bouquets $b \to b'$ whose root colors coincide is an arbitrary map $f: \underline{k} \to \underline{l}$ of their underlying sets. Otherwise, there is no map between $\frC$-bouquets. We denote the resulting category of $\frC$-bouquets by $\Bq(\frC)$.

The cardinality functor $\card{-}: \Bq(\frC) \to \Fin$ assigns to a bouquet $b: \underline{k} + 1 \to \frC$ its underlying set $\underline{k}$. The fiber of a map $b\to b'$ given by $f: \underline{k} \to \underline{l}$ over an element $y \in \underline{l}$ is a $\frC$-bouquet whose underlying set is $f^{-1}(y)$, the root color coincides with the color of $y$ and the colors of the elements are inherited from the colors of the elements of $\underline{k}$. It is easy to see that $\Bq(\frC)$ is an operadic category with $\frC$ as its set of connected components.
\end{example}

The category $\Bq(\frC)$ has the following important property. For each operadic category $\O$ with its set of connected components $\pi_0(\O)=\frC$, there is a canonical operadic \textit{Arity} functor $Ar_{\O}: \O\to \Bq(\frC)$ giving rise to the factorization
  \begin{equation*}
\begin{tikzcd}[ampersand replacement=\&,cramped]
	{\O} \&\& \Fin \\
	\& {\Bq(\frC)}
	\arrow["{\card{-}}", from=1-1, to=1-3]
	\arrow["Ar_{\O}"', from=1-1, to=2-2]
	\arrow["{\card{-}}"', from=2-2, to=1-3]
\end{tikzcd}
\end{equation*}
of the cardinality functor $\card{-} : \O \to \Fin$. 
We cite the construction of the $Ar_{\O}$ functor presented in \cite[Part I, Section 1]{opcat}. 
    Let the {\em source\/} $s(T)$ of $T \in \O$ be the set of fibers of the identity $id : T \to T$.  We define $Ar_{\O}(T) \in \Bq(\frC)$ as the
  bouquet $b: s(T)+1 \to \frC$, where $b$ associates to each fiber $U_c
  \in s(T)$ the corresponding connected component $c \in \frC$, and $b(1)
  := \pi_0(T)$.  The assignment
  $T \mapsto Ar_{\O}(T)$\/ extends into an operadic functor.



We denote by $\Bq$ the full subcategory of $\CatOp$ spanned by categories $\Bq(\frC)$, where $\frC \in \Set$. We observe that any operadic functor $F: \Bq(\frC) \to \Bq({\mathfrak D})$ is uniquely determined by an assignment of colors $f: \frC \to {\mathfrak D}$. Given an operadic functor $F: \O \to \opcat{R}$, there is a unique way to define the functor \[\Bq(F):~\Bq(\pi_0(\O)) \to \Bq(\pi_0(\opcat{R}))\] such that the diagram 
\begin{equation}\label{BqF}
\begin{tikzcd}[ampersand replacement=\&,cramped] 
	\O \& {\opcat{R}} \\
	{\Bq(\pi_0(\O))} \& {\Bq(\pi_0(\opcat{R}))}
	\arrow["F", from=1-1, to=1-2]
	\arrow["{Ar_{\O}}"', from=1-1, to=2-1]
	\arrow["{Ar_{\opcat{R}}}", from=1-2, to=2-2]
	\arrow["{\Bq (F)}"', dashed, from=2-1, to=2-2]
\end{tikzcd}
\end{equation}
commutes. The functor $F$ defines an assignment of colors $f: \pi_0(\O) \to \pi_0(\opcat{R})$ by $f(c) = \pi_0F(U_c)$, where $U_c$ is a local terminal object of $\O$ and so is $F(U_c)$, since $F$ preserves the chosen local terminals. This gives rise to the functor $\Bq (F)$. Therefore, the assignment $\O \mapsto \Bq(\pi_0(\O))$ is functorial, we denote it by
\[\AR: \CatOp \to \Bq.\] 

\begin{prop} \label{arity-inc-adj}
    The inclusion $i: \Bq \hookrightarrow \CatOp$ is the right adjoint to the $\AR: \CatOp \to \Bq$.
\end{prop}
\begin{proof}
    The components of the unit transformation are \[\eta_{\O} = Ar_{\O} : \O \to \Bq(\pi_0(\O)),\] for each $\O~\in~\CatOp$. The components of the counit transformation are \[\varepsilon_{\Bq(\frC)} = id_{\Bq(\frC)}:~\AR~\circ~i~(\Bq(\frC))~\to~\Bq(\frC).\] It is easy to see that both $\eta$ and $\varepsilon$ are natural transformations and satisfy the triangle identities.
\end{proof}

An $\O$-collection in a complete, cocomplete closed symmetric monoidal category $\V$ is a family $E = \{E(T)\}_{T \in \O}$ of objects of $\V$ indexed by the objects of the category $\O$. For an $\O\/$-collection $E$ and a morphism $f:T\to S$ in $\O$
let \[E(f) = \bigotimes_{i \in \card{S}} E({f^{-1}(i)}).\]

\begin{defn}
An $\O$-operad is  an $\O$-collection  $\oper{P} = \{\oper{P}(T)\}_{T
\in \O}$ in $\V$ together with units
  \[
  \eta_c: I \to \oper{P}(U_c),\ c \in \pi_0(\O),
  \]
and structure maps
  \[
  \mu_{\oper{P}}^f: \oper{P}(f) \otimes \oper{P}(S)\to \oper{P}(T),\ f:T\to S,
  \]
satisfying the axioms. 

\begin{enumerate}
    \item[(COMP)] Let $T \xrightarrow{f} S \xrightarrow{g} R$ be morphisms in $\opcat{O}$ and denote by $h$ their composition $g \circ f$. Then the following diagram of structure maps of $\oper{P}$ combined with the canonical isomorphisms of products in $\V$ commutes: 
\[\begin{tikzcd}[ampersand replacement=\&,cramped,column sep=tiny]
	{\displaystyle\bigotimes_{i \in \card{R}} \oper{P}(f_i) \otimes \oper{P}(g) \otimes \oper{P}(R)} \&\& {\oper{P}(h) \otimes \oper{P}(R)} \\
	{} \\
	{\displaystyle\bigotimes_{i \in \card{R}} \oper{P}(f_i) \otimes \oper{P}(S)} \& {\oper{P}(f) \otimes \oper{P}(S)} \& {\oper{P}(T)}
	\arrow["{\otimes_{i}\mu^{f_i} \otimes id}", from=1-1, to=1-3]
	\arrow["{id \otimes \mu^{g}}"', from=1-1, to=3-1]
	\arrow["{\mu^{h}}", from=1-3, to=3-3]
	\arrow["\cong"{marking, allow upside down}, draw=none, from=3-1, to=3-2]
	\arrow["{\mu^f}", from=3-2, to=3-3]
\end{tikzcd}\]
    \item[(UN1)]
    The composition
\[\begin{tikzcd}[ampersand replacement=\&,cramped,column sep=scriptsize]
	{\oper{P}(T)} \& {\displaystyle\bigotimes_{i \in \card{T}} I \otimes \oper{P}(T)} \& {\displaystyle\bigotimes_{i \in \card{T}} \oper{P}(U_{c_i}) \otimes \oper{P}(T)} \& {\oper{P}(id_T)\otimes \oper{P}(T)} \& {\oper{P}(T)}
	\arrow[from=1-1, to=1-2]
	\arrow[from=1-2, to=1-3]
	\arrow["{=}"{marking, allow upside down}, draw=none, from=1-3, to=1-4]
	\arrow["{\mu^{id_T}}", from=1-4, to=1-5]
\end{tikzcd}\]
    is the identity for each $T \in \opcat{O}$, as well as the identity is
    \item[(UN2)]
    the composition
\[\begin{tikzcd}[ampersand replacement=\&,cramped,column sep=scriptsize]
	{\oper{P}(T)} \& {\oper{P}(T) \otimes I} \& {\oper{P}(T)\otimes \oper{P}(U_c)} \& {\oper{P}(!_T)\otimes \oper{P}(U_c)} \& {\oper{P}(T),} \& {c := \pi_0(T).}
	\arrow[from=1-1, to=1-2]
	\arrow[from=1-2, to=1-3]
	\arrow["{=}"{marking, allow upside down}, draw=none, from=1-3, to=1-4]
	\arrow["{\mu^{!_T}}", from=1-4, to=1-5]
\end{tikzcd}\]
\end{enumerate}

A {\em morphism\/}  $\phi : \oper{P}' \to \oper{P}''$  of $\opcat{O}$-operads in $\V$ is a collection $\{\phi_T : \oper{P}'(T) \to \oper{P}''(T)\}_{T \in \O}$ of morphisms in $\V$ commuting with the structure maps. $\O$-operads in $\V$ form a category $\Opc{\O}{\V}$.
\end{defn}

\begin{example}\label{ex-term-op}
    The category of $\O$-operads in $\Set$ has a terminal object, namely the operad $\termop_{\opcat{O}} \in \Opc{\opcat{O}}{\Set}$, where $\termop_{\opcat{O}}(T) = \{T\}, $ for $ T \in \opcat{O}$.
\end{example}

\begin{prop}[{\cite[Prop.~3.1.]{Batanin-EH}}] \label{Fin-op-is-op}
The category of classical operads in $\V$ in the spirit of May \cite{may} is isomorphic to the category of {\Fin}-operads in $\V$ in the sense of \cite{opcat}.
\end{prop}
\begin{proof}
We recall the correspondence and refer the reader to \cite{Batanin-EH} for more details. Suppose $\oper{P} \in \Opc{\Fin}{\V}$, we define the structure of a symmetric operad on $\oper{P}$. The structure map \[\gamma: \oper{P}(k) \otimes \oper{P}(n_1) \otimes \cdots \otimes \oper{P}(n_k) \to \oper{P}(n_1 + \cdots + n_k)\] is given by $\mu^{\nu}$, where $\nu: n_1 + \cdots + n_k \to k$ is an order-preserving morphism such that $\nu(n_i) = i$. The right action of $\pi \in \Sigma_n$ on $\oper{P}(n)$is given as the composite 
\[\begin{tikzcd}[ampersand replacement=\&]
	{\oper{P}(n)} \& {I^{\otimes n} \otimes \oper{P}(n)} \& {\oper{P}(1)^{\otimes n} \otimes \oper{P}(n) } \& {\oper{P}(n)}
	\arrow["\cong", from=1-1, to=1-2]
	\arrow["{\eta^n \otimes id}", from=1-2, to=1-3]
	\arrow["{\mu^{\pi}}", from=1-3, to=1-4]
\end{tikzcd}.\]
In case $\V = \Set$, let $u \in \oper{P}(1)$ be the image of the unit morphism $\eta: Pt \to \oper{P}(1)$. Then this translates to \[\alpha \cdot \pi := \mu^{\pi}((u,\ldots,u), \alpha),\]  for $\alpha \in \oper{P}(n)$ and $\pi \in \Sigma_n$.

In the other direction, suppose $\oper{Q}$ is a symmetric operad in $\V$. We define the structure of a $\Fin$ operad on $\oper{Q}(n)$ as follows. To define the composition $\mu^{\sigma}$ along a morphism $\sigma: n \to m$ in $\Fin$, we recall that every such morphism has a unique decomposition
\[\begin{tikzcd}[ampersand replacement=\&]
	n \&\& k \\
	\& {n'}
	\arrow["\sigma", from=1-1, to=1-3]
	\arrow["{\pi(\sigma)}"', from=1-1, to=2-2]
	\arrow["{\nu(\sigma)}"', from=2-2, to=1-3]
\end{tikzcd}\]
into a permutation $\pi(\sigma)$ and an order-preserving $\nu(\sigma)$ such that the order of fibers is preserved. We use this factorization to define $\mu^{\sigma}((\alpha_1,\ldots, \alpha_k),\beta) := \gamma(\beta, \alpha_1,\ldots, \alpha_k) \cdot \pi(\sigma)$.
\end{proof} 

The following generalization of Proposition \ref{Fin-op-is-op} holds by the same arguments.

\begin{prop} \label{col-op-bq-op}
The category of operads over the category $\Bq(\frC)$ of $\frC$-bouquets 
is isomorphic to the category of $\frC$-colored symmetric operads.
\end{prop}

Observe that an operadic functor $F: \O \to \opcat{R}$ induces the restriction $F^*:\Opc{\opcat{R}}{\V} \to \Opc{\O}{\V}$, where $F^*(\oper{P})(T) = \oper{P}(F(T))$ and $\mu_{F^*(\oper{P})}^f = \mu_{\oper{P}}^{Ff}$. We introduce an important class of operadic functors such that the restriction $F^*$ has a left adjoint $F_!$. We say that an operadic functor $F: \opcat{O} \to \oper{P}$ is \textit{a discrete operadic fibration} if
\begin{enumerate}
    \item $F$ induces an epimorphism $\pi_{0}(\O) \twoheadrightarrow \pi_0(\opcat{R})$;
    \item for any morphism $f: T \to S$ in $\opcat{R}$ and $t_i, s \in \opcat{O}$, where $i \in \card{S}$ such that 
    \[F(s) = S \text{~~and~~} F(t_i) = f^{-1}(i),\] there exists a unique $\sigma: t \to s$ in $\opcat{O}$ such that \[F(\sigma) = f  \text{~~and~~} t_i = \sigma^{-1}(i).\] 
\end{enumerate}

Given a discrete operadic fibration $F: \O \to \opcat{R}$ and an operad $\oper{P} \in \Opc{\O}{\V}$, the collection $F_! (\oper{P})$
\[F_! (\oper{P}) (T) = \{\coprod_{F(t) = T} \oper{P}(t)\},\] for $T \in \opcat{R}$, has a natural $\opcat{R}$-operad structure \cite[Prop.~2.3.]{opcat}, which defines the left adjoint $F_!$ to the restriction $F^*$ \cite[Prop.~2.4.]{opcat}.

\begin{notat}
    From now on, the ambient category $\V$ will be the category of $\Set$ of sets, and we will omit all indices referring to a specific ambient category.
\end{notat}

Another class of functors, for which the induced restriction has a left adjoint, is a class of arity functors $Ar_{\O}: \O \to \Bq(\pi_0(\O))$. To construct  \[Ar^{\O}_!: \Opc{\O}{} \to \Opc{\Bq(\pi_0(\O))}{},\] given a $\O$-operad $\oper{Q}$, we define a $\Bq(\pi_0(\O))$-collection by
\begin{equation*} \label{eq-coll}
    E_{\oper{Q}}(T) := \coprod_{Ar_{\O}(t) = T} \oper{Q}(t)
\end{equation*} for each $T \in \Bq(\pi_0(\O))$. Denote $\oper{F}_{\oper{Q}}$ the free unital colored operad generated by the collection~$E_{\oper{Q}}$. We identify the new adjoint units $u'_w \in \oper{F}_{\oper{Q}} (\overset{w}{\underset{w}{\mid}})$ with the units $u_w \in \oper{Q}(U_w)$, for each $w$ connected component of $\opcat{O}$ and $U_w$ its local terminal object. Then, we take the quotient of $\oper{F}_{\oper{Q}}$ by the operadic ideal generated by pairs of the form
\begin{equation} \label{eq-rel}
    \mu_{\oper{F}_{\oper{Q}}} ^{Ar_{\O}(f)} ((y_1,...,y_n),x) \sim z,
\end{equation} where $Y_1, \ldots, Y_n\rhd Z \xrightarrow{f} X$ is a morphism in $\O$ and $x \in \oper{Q}(X), y_i \in \oper{Q}(Y_i), z \in \oper{Q}(Z)$ such that \[\mu_{\oper{Q}} ^{f} ((y_1,...,y_n),x) = z\] holds in $\oper{Q}$. We define 
\begin{equation}\label{constr-ar!}
    Ar^{\O}_!(\oper{Q}) := \faktor{\oper{F}_{\oper{Q}}}{\sim}.
\end{equation}
It is straightforward to show that the construction above is functorial.

\begin{prop}
    The functor $Ar^{\O}_{!}: \Opc{\O}{} \to \Opc{\Bq(\pi_0(\O))}{}$ defined above is the left adjoint to the restriction functor $Ar^*_{\O}:  \Opc{\Bq(\pi_0(\O))}{} \to \Opc{\O}{}$.
\end{prop}
\begin{proof}
    We show that there is a natural bijection of sets \[\Opc{\Bq(\pi_0(\O))}{}(Ar^{\O}_!(\oper{Q}), \oper{P}) \cong \Opc{\O}{}(Q, Ar_{\O}^*(\oper{P})),\] for each $\oper{Q} \in \Opc{\O}{}$ and $\oper{P} \in \Opc{\Bq(\pi_0(\O))}{}$. Let $\varphi: \oper{Q} \longrightarrow Ar_{\O}^*(\oper{P})$ be a morphism of $\opcat{O}$-operads. It consists of components \[\varphi_T: \oper{Q}(T) \longrightarrow Ar_{\O}^*(\oper{P})(T) = \oper{P}(Ar_{\O}(T)),\] for each $T \in \opcat{O}$, that assemble to \[\coprod_{Ar_{\O}(T) = t} \varphi_T : \coprod_{Ar_{\O}(T) = t} \oper{Q}(T) \to \oper{P}(t),\] for each $t \in \Bq(\pi_0(\O))$. This extends to a morphism from the free operad $\tilde{\varphi}: \oper{F}_{\oper{Q}} \longrightarrow \oper{P}.$ We need to show that the equivalence relation (\ref{eq-rel}) is in the kernel of $\tilde{\varphi}$.

    Suppose $y_1,...,y_n,x$ are as described in (\ref{eq-rel}). Then
    \[\arraycolsep=2pt\def\arraystretch{1.6}
    \begin{array}{rcll}
        \tilde{\varphi}(\mu_{\oper{F}_{\oper{Q}}}^{Ar_{\opcat{O}}f}((y_1,...,y_n),x)) & = & \mu_{\oper{P}}^{Ar_{\opcat{O}}f}((\tilde{\varphi}(y_1),...,\tilde{\varphi}(y_n)),\tilde{\varphi}(x)) & \text{(since $\tilde{\varphi}$ is a morphism of operads)} \\
         & = & \mu_{\oper{P}}^{Ar_{\opcat{O}}f}(({\varphi}(y_1),...,{\varphi}(y_n)),{\varphi}(x))
         & \text{(by def. of $\tilde{\varphi}$ on generators of $\oper{F}_{\oper{Q}}$)}\\
         & = & \mu_{Ar_{\O}^*\oper{P}}^{f}(({\varphi}(y_1),...,{\varphi}(y_n)),{\varphi}(x)) & \text{(by def. of restriction $Ar^*_{\opcat{O}}$)}\\
         & = & \varphi(\mu_{\oper{Q}}^{f}((y_1,...,y_n),x)) & \text{(since $\varphi$ is a morphism of operads)} \\
         & = & \varphi(z) = \tilde{\varphi}(z).
    \end{array}
    \]
    This verifies that $\tilde{\varphi}$ factors through the morphism $\varphi^{\musSharp}: Ar^{\O}_!(\oper{Q}) \to \oper{P}$ defined by the assignment $\varphi^{\musSharp}([x]) = \tilde{\varphi}(x)$.

    In the opposite direction, let $\psi: Ar^{\opcat{O}}_!(\oper{Q}) \to \oper{P}$ be a morphism of $\Bq(\pi_0(\O))$-operads. Let $x \in \oper{Q}(X)$; then the equivalence class $[x]$ under the relation (\ref{eq-rel}) is an element in $Ar_!^{\O}(\oper{Q})(Ar_{\O}(X))$, and $\psi([x])$ is an element in $\oper{P}(Ar_{\O}(X))$. We define a morphism $\psi^{\musFlat}: \oper{Q} \to Ar^*_{\O}(\oper{P})$ of $\O$-operads by the assignment $\psi^{\musFlat}(x) = \psi([x])$.

    We need to show that $\psi^{\musFlat}$ is a morphism of operads. Assume $Y_1, \ldots, Y_n\rhd Z \xrightarrow{f} X$ is a morphism in $\O$ and $x \in \oper{Q}(X), y_i \in \oper{Q}(Y_i), z \in \oper{Q}(Z)$ such that $\mu_{\oper{Q}} ^{f} ((y_1,...,y_n),x) = z$ holds in $\oper{Q}$.
    \begin{multline*} 
        \mu^f_{Ar^*_{\O}(\oper{P})} ((\psi^{\musFlat}(y_1), \ldots,\psi^{\musFlat}(y_n)), \psi^{\musFlat}(x)) = \\
        \arraycolsep=2pt\def\arraystretch{1.6}
        \begin{array}{rll}
            = & \mu^{Arf}_{\oper{P}} ((\psi[y_1], \ldots,\psi[y_n]), \psi[x]) & \text{(by def. of $Ar^*_{\O}$ and $\psi^{\musFlat}$)} \\
            = & \psi(\mu^{Arf}_{Ar^{\O}_!(\oper{Q})} (([y_1], \ldots, [y_n]), [x])) & \text{(since $\psi$ is a morphism of operads)} \\
            = & \psi([z]) & \text{(by definition of composition in $Ar^{\O}_!(\oper{Q})$)} \\
            = & \psi^{\musFlat}(\mu_{\oper{Q}} ^{f} ((y_1,...,y_n),x)).
        \end{array}
    \end{multline*}
It is straightforward to show that the assignments above are inverse to each other and that the bijection is natural.
\end{proof}

In case $Ar_{\O}: \O \to \Bq(\pi_0(\O))$ is a discrete operadic fibration, all free compositions are equivalent to some element of the operad $\oper{Q}$. Hence, the components of $Ar^{\opcat{O}}_!(Q)$ are just coproducts of fibers, and the structure of the $\Bq(\pi_0(\O))$-operad is the natural one induced by the discrete operadic fibration described in \cite[Prop.~2.3.]{opcat}.  

Let $\opcat{O}$ be an operadic category and an $\oper{P} \in \Opc{\O}{}$. The \textit{operadic Grothendieck construction} \cite[Prop.~2.5.]{opcat} is the category $\displaystyle\int_{\opcat{O}} \oper{P}$ whose objects are $t \in \oper{P}(T)$ for some $T \in \opcat{O}$. A morphism $\sigma : t \to s$ from $t \in \oper{P}(T)$ to $s \in \oper{P}(S)$ is a pair $(\varepsilon,f)$
consisting of a morphism $f : T \to S$ in $\opcat{O}$ and a tuple
$\varepsilon \in \prod_{i \in |S|} \oper{P}(f^{-1}(i))$, such that
$
\mu_{\oper{P}}^f(\varepsilon,s) = t,
$
where $\mu_{\oper{P}}$ is the structure map of the operad $\oper{P}$.  Compositions of
morphisms are defined in the obvious manner.  

%
%

We use $\Op$ to denote \textit{the category of operads} in $\Set$. The objects of $\Op$ are pairs $(\oper{P} \in \Opc{\O}{})$, where $\O \in \CatOp$. A morphism $(\oper{P} \in \Opc{\O}{}) \to (\oper{Q} \in \Opc{\opcat{R}}{})$ consists of a pair $F: \O \to \opcat{R}$ in $\CatOp$ and $f: \oper{P} \to F^*(\oper{Q})$ in $\Opc{\O}{}$. Denote by $\SOp$ the full subcategory of symmetric unital colored operads of $\Op$, thus objects of $\SOp$ are pairs $(\oper{B} \in \Opc{\Bq(\frC)}{})$, for some set of colors $\frC \in \Set$.

We prove that $\SOp$ is a reflective subcategory of $\Op$, i.e. construct a left adjoint $\Arleft: \Op \to \SOp$ to the inclusion $inc: \SOp \to \Op$. We define the action on objects to be $\Arleft(\oper{P} \in \Opc{\O}{}) := Ar^{\O}_{!}(\oper{P})$. Suppose that $(F, f)$ is a morphism $(\oper{P} \in \Opc{\O}{}) \to (\oper{Q} \in \Opc{\opcat{R}}{})$. 
The morphism \[\Arleft(F, f): (Ar^{\O}_{!}(\oper{P}) \in \Opc{\Bq(\pi_0(\opcat{O}))}{}) \to (Ar^{\opcat{R}}_{!}(\oper{Q}) \in \Opc{\Bq(\pi_0(\opcat{R}))}{})\] consists of a functor $\Bq(F): \Bq(\pi_0(\opcat{O})) \to \Bq(\pi_0(\opcat{R}))$ and a morphism \[x:Ar^{\O}_{!}(\oper{P}) \to \Bq(F)^*\circ Ar^{\opcat{R}}_{!}(\oper{Q})\] in $\Opc{\Bq(\pi_0(\opcat{O}))}{}$. Since $Ar^{\O}_{!} $ is the left adjoint to $Ar^{*}_{\O}$, to specify $x$, it is enough to specify \[x^{\musFlat}: \oper{P} \to Ar^{*}_{\O}\circ \Bq(F)^*\circ Ar^{\opcat{R}}_{!}(\oper{Q})\] in $\O$. However, since $\Bq(F)$ is such that (\ref{BqF}) commutes, the equality $Ar^{*}_{\O}\circ \Bq(F)^*= F^*\circ Ar^*_{\opcat{R}}$ holds. We define $\underline{x}$ to be the composite
\[x^{\musFlat}:\oper{P} \xrightarrow{~f~} F^*(\oper{Q})
                   \xrightarrow{~F^*(\eta_{\oper{Q}}^{\opcat{R}})~}
                   F^* \circ Ar^*_{\opcat{R}} \circ Ar_!^{\opcat{R}} (\oper{Q}),\]
where $\eta^{\opcat{R}}$ is the unit of the adjunction $Ar^{\opcat{R}}_{!} \dashv Ar^{*}_{\opcat{R}}$. It is straightforward to show that $\Arleft$ is a functor.

\begin{prop} \label{arleft-inc-adj}
    There is an adjunction 
    \[\xymatrix@R=1em@C=1.5em{{\Op} \ar@/^0.9pc/[r]<1ex>^{\Arleft} & \SOp\ar@/^0.9pc/[l]<1ex>^{inc} \ar@{}[l]|{\upvdash}}.\] 
\end{prop}
\begin{proof}
The components of the unit transformation \[\eta_{\oper{P}}: \oper{P} \to inc \circ \Arleft(\oper{P}), \quad \oper{P} \in \Opc{\O}{}\] are pairs $Ar_{\opcat{O}}: \opcat{O} \to \Bq(\pi_0(\O))$ and $\eta^{\opcat{O}}_{\oper{P}}: \oper{P} \to Ar_{\opcat{O}}^* Ar^{\opcat{O}}_! (\oper{P})$, for $\oper{P} \in \Opc{\O}{}$. For the counit transformation, we observe that $\Arleft \circ inc$ is an identity functor since $Ar_{\Bq(\frC)}$ is an identity for any bouquet category $\Bq(\frC)$. We define the counit to be the identity transformation. It is straightforward to verify the triangle identities.
\end{proof}

\begin{rem}
    Consider the functor $\text{Op}: \CatOp^{op} \to \text{CAT}$ that assigns to a category $\opcat{O}$ the category $\Opc{\O}{}$, and whose action on operadic functors is given by restriction. The category $\Op$ is then the Grothendieck fibration associated to $\text{Op}$. Similarly, the category $\SOp$ is the Grothendieck fibration associated to the restriction of $\text{Op}$ to the category of bouquet operadic categories $\Bq$. Readers familiar with base changes for adjunctions (see, for example, \cite{grothendieck}) may recognize in this setup the construction of the base change for the adjunction $\AR \dashv i$ in Proposition \ref{arity-inc-adj}.
    
     To apply this base change more generally, however, one would need the functor $\text{Op}$ to be a bifibration. Proposition \ref{arleft-inc-adj} suggests that $\text{Op}$ is indeed likely to be a bifibration, and therefore that each restriction $F^*$ admits a left adjoint $F_!$. A full proof of this assertion, however, lies beyond the scope of the present paper.
\end{rem}

\begin{prop} \label{CatOp-Op-adj}
    There is an adjunction of categories 
    \[\xymatrix@R=1em@C=1.5em{{\CatOp~~}\ar@/^0.9pc/[r]<1ex>^{\termop} & \Op\ar@/^0.9pc/[l]<1ex>^G \ar@{}[l]|{\upvdash}}.\] 
    The right adjoint  is the operadic Grothendieck construction
    \[\begin{array}{rccc}
         G: & \Op & \to  &\CatOp  \\
         &  \oper{P} \in \Opc{\O}{} & \longmapsto & \int_{\O} \oper{P} 
    \end{array}\]
    and the left adjoint is the  terminal operad over an operadic category
    \[
    \begin{array}{rccc}
         \termop: & \CatOp  & \to  & \Op  \\
         &  \O & \longmapsto & \termop_{\O} \in \Opc{\O}{}
    \end{array}
    \]
\end{prop}
\begin{proof}
    The collection of isomorphisms \[\eta_{\O}: \O \xrightarrow{~\sim~} \int_{\O} \termop_{\O},\] for each $\O \in \CatOp$, defines the unit transformation. A component of the counit transformation,
    \[\varepsilon_{\oper{P}}: \termop \circ G(\oper{P} \in \Opc{\O}{}) \to \oper{P} \in \Opc{\O}{},\] for $\oper{P} \in \Opc{\O}{}$, is given by a projection $\pi: \int_{\O} \oper{P} \to \O$ and a morphism $p: \termop_{\int_{\opcat{O}} \oper{P}} \to \pi^* \oper{P}$, where
    \[p_x: \termop_{\int_{\opcat{O}} \oper{P}}(x) \to \pi^* \oper{P}(x) = \oper{P}(T), \quad \text{for } x \in \oper{P}(T)\] is the inclusion $\{x\} \hookrightarrow \oper{P}(T)$.
\end{proof}

We, therefore, have a chain of adjunctions
\[\xymatrix@R=1em@C=1.5em{{\CatOp~~}\ar@/^0.9pc/[r]<1ex>^{\termop} & \Op\ar@/^0.9pc/[l]<1ex>^G \ar@{}[l]|{\upvdash} \ar@/^0.9pc/[r]<1ex>^{\Arleft} & \SOp\ar@/^0.9pc/[l]<1ex>^{inc} \ar@{}[l]|{\upvdash}}.\] 

We denote $\A = \Arleft \circ \termop$ and $\I = G \circ inc$.

\begin{obs} \label{obs-AI-op}
    Given a $\frC$-colored operad $\oper{P}$, the composite $\A \I(\oper{P})$ is canonically isomorphic to~$\oper{P}$. Indeed, $Ar: \int_{\Bq(\frC)} \oper{P} \to \Bq(\frC)$ is a discrete operadic fibration by definition of the operadic Grothendieck construction; therefore, the components of $\A \I (\oper{P})$ is given only by coproducts of fibers of $Ar$. If we apply this process to the terminal operad, we reconstruct the operad $\oper{P}$.
\end{obs}

We formulate the main result of this section.

\begin{thm} \label{thm-adj-AI}
    The category $\opcat{SO}$ of symmetric $\frC$-colored operads in $\Set$ is a reflective subcategory of $\opcat{CatOp}$ of operadic categories.
\end{thm}

%% file: wreath_prod.tex
\begin{defn} \label{defn-prewr}
Suppose $\opcat{O}, \opcat{R}$ are operadic categories. We define their \emph{operadic wreath product} $\OwR$ to be the category
\begin{itemize}
    \item with objects the tuples $(X; Y_1, \ldots, Y_n)$, where $X \in \opcat{O}$, $\card{O} = \underline{n}$, and $Y_1, \ldots Y_n \in \opcat{R}$ such that they belong to the same connected component of $\opcat{R}$;
    \item with morphisms $(\varphi, \Phi): (X; Y_1,\ldots, Y_n) \to (Z; W_1, \ldots W_k)$, which consist of a morphism $\varphi: X \to Z$ in $\opcat{O}$, and a family $\Phi = \{\varphi_{ij}: Y_i \to W_j \mid \card{\varphi}(i) = j\}$ of morphisms in $\opcat{R}$. 
\end{itemize}
\end{defn}

We give $\OwR$ the structure of an operadic category in the following way. We define the cardinality functor on objects \[\card{(X; Y_1,\ldots, Y_n)}~:=~\coprod_{i \in \underline{n}} \card{Y_i},\] and on morphisms \[\card{(\varphi, \Phi)}~:=~\coprod_{\substack{i \in \underline{n} \\ \varphi(i) = j}} \card{\varphi_{ij}}.\]
Similarly, the fibers of $\opcat{O}$ and $\opcat{R}$ induce pointwise fibers in $\OwR$, the $i$-th fiber is the object \[(\varphi, \Phi)^{-1}(i)= (\varphi^{-1}(t); ~ \big(\varphi_{st}^{-1}(i)\big)_{s \in \card{\varphi}^{-1}(t)} ~).\] The following observation implies that  this definition indeed produces an object of $\OwR$, i.e. that the pointwise fibers in $\opcat{R}$ lie in the same connected component.

\begin{obs} \label{obs-fib}
Suppose we are given two morphisms $f: Y_1 \to W$ and $g: Y_2 \to W$. Then, by (AX2) and (AX3), the commutative triangle 
\[\begin{tikzcd}[ampersand replacement=\&,cramped]
	{Y_1} \&\& W \\
	\& W
	\arrow["f", from=1-1, to=1-3]
	\arrow["f"', from=1-1, to=2-2]
	\arrow["id", from=1-3, to=2-2]
\end{tikzcd}\]
induces a morphism $f_i: f^{-1}(i) \to U$, where $i \in \card{w}$ and $U$ is some chosen local terminal object. Similarly, there exists a morphism $g_i: g^{-1}(i) \to U$. Therefore, $f^{-1}(i)$ and $g^{-1}(i)$ lie in the same connected component.
\end{obs}

We chose the family of local terminal objects to be 
\[
\{(U,V) \mid U, \text{ resp. } V,\text{ is a chosen local terminal of } \opcat{O}, \text{ resp. } \opcat{R}\}.
\]

\begin{prop}
    The category $\OwR$ is an operadic category.
\end{prop}
\begin{proof}
    This is a straightforward verification done in \cite{masters}. 
\end{proof}

It is clear that the wreath product of operadic categories is not commutative but is associative; the natural isomorphism $\wprod{(\OwR)}{\opcat{Q}} \cong \wprod{\opcat{O}}({\wprod{\opcat{R}}{\opcat{Q}}})$ is given by simple rebracketing. Moreover, given an operadic category $\opcat{O}$ and an operadic category $\mathbf{1}$ consisting of one object $\ast$ of cardinality $\underline{1}$ and an identity morphism, there exist obvious natural isomorphisms $\wprod{\mathbf{1}}{\opcat{O}} \cong \opcat{O}$ and $\wprod{\O}{\mathbf{1}} \cong \O$ such that the known pentagon and triangle identities are satisfied. This allows us to conclude the following. 

\begin{prop}
    $(\CatOp, \wr, \mathbf{1})$ is a monoidal category.
\end{prop}

The category  $\Omega_k$ of $k$-trees, $k\geq 0$, \cite[Sec.~3, Ex.~8]{BATANIN199839} is an operadic category, as shown in \cite[Sec.~1.1.]{opcat}. We observe that it can be constructed iteratively using the wreath product of operadic categories. We recall that $\Omega_1 := \Dalg$.

\begin{prop} \label{prop-trees}
   Let $k \geq 1$, then $\wprod{\Omega_{k}}{\Omega_1} \cong \Omega_{k+1}$.    
\end{prop}
\begin{proof}
    Given $((n_k \xrightarrow{t_{k-1}}  \ccdots  \xrightarrow{t_1} n_1); (a_1),  \ccdots , (a_{n_k}))$ in $\wprod{\Omega_{k}}{\Omega_1}$, the isomorphism $\wprod{\Omega_{k}}{\Omega_1}~\xrightarrow{\sim}~\Omega_{k+1}$ acts by constructing a tree \[\displaystyle\coprod_{i=1}^{n_k} (a_i) \xrightarrow{t_k} n_k \xrightarrow{t_{k-1}}  \ccdots  \xrightarrow{t_1} n_1,\] where $t_k(a_i) = i$.
%
%
%
%
%
%
%
%
%
%
%
%
%
\end{proof}

We now aim to show that the functor $\A: (\CatOp, \wr, \mathbf{1}) \to (\SOp, \BV, \star)$ is a monoidal functor, i.e. that there is an isomorphism
\begin{equation*}
\iota: \A(\mathbf{1}) \xrightarrow{\sim} \star    
\end{equation*}
and natural isomorphisms, for all $\opcat{O}, \opcat{R} \in \CatOp$,
\begin{equation*}\gamma_{\opcat{O}, \opcat{R}}: \A(\opcat{O} \wr \opcat{R}) \xrightarrow{\sim} \A(\opcat{O}) \BV \A(\opcat{R}),
\end{equation*}
satisfying the classic conditions.
Here, the monoidal structure on the category $\SOp$ is given by the Boardman-Vogt tensor product of operads  $\BV$, and the initial unital operad $\star$ \cite[Th.~2.22.]{weiss}. We first recall the definition of the BV product in terms of generators and relations, as presented, for example, in \cite[Def.~2.21.]{weiss}.

\begin{defn} \label{def-bv}
    Let $\oper{P}$ be a symmetric $\frC$-colored operad  and $\oper{Q}$ be a symmetric $\frD$-colored operad. Their \textit{Boardman-Vogt tensor product} is the symmetric operad \hbox{$\oper{P} \BV \oper{Q}$} with a set of colors $\frC \times \frD$. The operad \hbox{$\oper{P} \BV \oper{Q}$} is generated by two families of generators:
    \begin{itemize}
        \item generators of the type $x \otimes d \in (\oper{P} \BV \oper{Q})\begin{pmatrix}
        (c_1,d), \ldots, (c_n,d) \\
        (c,d)
    \end{pmatrix}$, for each $x \in \oper{P}\begin{pmatrix}
        c_1 \ldots c_n \\
        c
    \end{pmatrix}$ and each color $d \in \frD$;
    \item generators of the type $c \otimes y\in (\oper{P} \BV \oper{Q})\begin{pmatrix}
        (c,d_1), \ldots, (c,d_m) \\
        (c,d)
    \end{pmatrix}$, for each color $c \in \frC$ and each $y \in \oper{Q}\begin{pmatrix}
        d_1 \ldots d_m \\
        d
    \end{pmatrix}$,
    \end{itemize}
    so that the following relations hold:
    \begin{itemize}
        \item $\gamma_{\oper{P} \BV \oper{Q}} (x\otimes d, x_1 \otimes d, \ldots,x_n\otimes d) = \gamma_{\oper{P}}(x, x_1,\ldots, x_n) \otimes d,$ for composable \hbox{$x, x_1, \ldots, x_n \in \oper{P}$} and $d \in \frD$;

        \item $(x \cdot \sigma) \otimes d = (x \otimes d) \cdot \sigma,$ for $x \in \oper{P}$, $d \in \frD$ and an appropriate permutation $\sigma$;
    \end{itemize}
    This ensures that  for any color $d \in \frD$ the inclusion $-\otimes d: \oper{P} \hookrightarrow \oper{P} \BV \oper{Q}$ given by $x \mapsto x \otimes d$ is a morphism of operads. Similar relations ensure that, for any color $c \in \frC$, the inclusion $c \otimes -: \oper{Q} \to \oper{P} \BV \oper{Q}$ is a morphism of operads.

    Lastly, the \textit{interchange} relation must hold, i.e., for any $x \in \oper{P}\begin{pmatrix}
        c_1 \ldots c_n \\
        c
    \end{pmatrix}$ 
    and 
    $y \in \oper{Q}\begin{pmatrix}
        d_1 \ldots d_m \\
        d
    \end{pmatrix}$,
    \[\gamma_{\oper{P} \BV \oper{Q}} (x \otimes d, c_1 \otimes y, \ldots, c_n \otimes y) = \gamma_{\oper{P} \BV \oper{Q}} (c \otimes y, x \otimes d_1, \ldots, x \otimes d_n) \cdot \textit{shuffle},\] where \textit{shuffle} is the permutation, the role of which we illustrate below. Consider the expressions
    \begin{center}
        \includegraphics[width=0.9\linewidth]{interchange_prelim.tex}.
    \end{center}
    The compositions on the left-hand side and the right-hand side cannot be identified since their domains differ. For this reason, we apply the \textit{shuffle} permutation to the composition on the right-hand side, which reorders the colors from reverse-lexicographical order to lexicographical order.
\end{defn}

We first examine the object $\A(\mathbf{1}) := Ar^{\mathbf{1}}_!\termop_{\mathbf{1}}$. By definition, recall Ex. \ref{ex-term-op}, $\termop_{\mathbf{1}}(\ast) = \{\ast\}$. Since the category $\mathbf{1}$ has only one morphism $\ast \xrightarrow{id} \ast$, the only structure morphism is the trivial composition given by $\mu^{id}_{\termop_{\mathbf{1}}}(\ast, \ast) = \ast$. Next, we apply the functor $Ar^{\mathbf{1}}_!$ as constructed in (\ref{constr-ar!}). 

Following the notation of (\ref{constr-ar!}), we start with the generating $\Fin$-collection \[E_{\termop_{\mathbf{1}}}(T) = \begin{cases}
    \{\ast\}, & T = \underline{1} \\
    \emptyset, & \text{otherwise.}
\end{cases}\] The free $\Fin$-operad $\oper{F}_{\termop_{\mathbf{1}}}$ is the free one-colored unital symmetric operad on one generator in arity~1. After identifying the freely adjoint unit with the element $\ast \in \oper{F}_{\termop_{\mathbf{1}}} (\underline{1})$ we are left with the operad isomorphic to the initial operad $\star$. We denote the unique isomorphism by 
$\iota: \A(\mathbf{1}) \xrightarrow{\sim} \star.$

Our next goal is to define the natural transformation $\gamma: \A(- \wr -) \Rightarrow \A(-) \BV \A(-)$. Fix any such $\opcat{O}, \opcat{R}$. Since $\SOp$ is a reflective subcategory of $\Op$, given any morphism \[(F_{\opcat{O}, \opcat{R}}, f_{\opcat{O}, \opcat{R}}): \termop_{\OwR} \to \A(\opcat{O}) \BV \A(\opcat{R}),\] there exists a unique morphism \[\gamma_{\opcat{O}, \opcat{R}}: \A(\OwR) \to \A(\opcat{O}) \BV \A(\opcat{R})\] such that the following triangle commutes.
\[\begin{tikzcd}[ampersand replacement=\&,cramped]
	{\termop_{\OwR}} \& {\A(\opcat{O}) \BV \A(\opcat{R})} \\
	{\A(\OwR)}
	\arrow["{(F_{\opcat{O}, \opcat{R}}, f_{\opcat{O}, \opcat{R}})}", from=1-1, to=1-2]
	\arrow["\eta"', from=1-1, to=2-1]
	\arrow["{\gamma_{\opcat{O}, \opcat{R}}}"', dotted, from=2-1, to=1-2]
\end{tikzcd}\]
We put $F_{\opcat{O}, \opcat{R}}:= Ar_{\OwR}: \OwR \to \Bq(\pi_0{\opcat{O}} \times \pi_0{\opcat{R}})$ and want to define \[f_{\opcat{O}, \opcat{R}}: \termop_{\OwR} \to Ar_{\OwR}^{\ast}(\A(\opcat{O}) \BV \A(\opcat{R})).\] Given $(X; Y_1, \ldots, Y_n) \in \termop_{\OwR}(X; Y_1, \ldots, Y_n)$, where 
$Ar_{\opcat{O}}(X)= (c_1, \ldots, c_n;c)\text{, and } Ar_{\opcat{R}}(Y_i) = (d_1^i, \ldots, d_{\card{y_i}}^i)$,
%
%
we observe that \[Ar_{\OwR}(X; Y_1, \ldots, Y_n)= \big(
(c_1, d_1^1), (c_1, d_2^1), \ldots, (c_n, d_{\card{Y_n}-1}^n), (c_n, d_{\card{Y_n}}^n);(c,d)
\big)\]
%
Moreover, $X \otimes d$ is an element of $\A(\opcat{O}) \BV \A(\opcat{R})\begin{pmatrix}
        (c_1,d), \ldots, (c_n,d) \\
        (c,d)
    \end{pmatrix}$ and similarly, $c_i \otimes Y_i$ is an element of $\A(\opcat{O}) \BV \A(\opcat{R})\begin{pmatrix}
        (c_i,d^i_1), \ldots, (c_i,d^i_{\card{y_i}}) \\
        (c_i,d)
    \end{pmatrix}$.
There is an obvious morphism of $\pi_0\opcat{O} \times \pi_0\opcat{R}$-colored bouquets.
\[c_i \otimes Ar_{\opcat{R}}(Y_i) \rhd Ar_{\OwR}(X; Y_1, \ldots, Y_n) \xrightarrow{\nu} Ar_{\opcat{O}}(X) \otimes d, \] along which the composition
\[\mu_{BV}^{\nu}\left(\left(c_1 \otimes Y_1, \ldots, c_n \otimes Y_n\right), X \otimes d\right)\]
produces an element of $\A(\opcat{O}) \BV \A(\opcat{R}) (Ar_{\OwR}(X; Y_1, \ldots, Y_n)) = Ar^{\ast}_{\OwR}(\A(\opcat{O}) \BV \A(\opcat{R}))(X; Y_1, \ldots, Y_n)$. 
We put 
\begin{equation} \label{morph-def}
    f_{\opcat{O}, \opcat{R}}(X; Y_1, \ldots, Y_n) := \mu_{BV}^{\nu}\left(\left(c_1 \otimes Y_1, \ldots, c_n \otimes Y_n\right), X \otimes d\right).
\end{equation}

\begin{notat}
    From now on, for finite linearly ordered sets $\underline{n}$ and $\underline{p}_1, \ldots, \underline{p}_n$ we denote by $\nu: \coprod_{i \in \underline{n}} \underline{p}_i \to \underline{n}$ the order-preserving map $\nu(\underline{p}_i) = i$. This notation extends in the obvious way to maps of bouquets with an underlying order-preserving map.
\end{notat}

\begin{prop}
    The assignment (\ref{morph-def}) is a morphism of operads.
\end{prop}
\begin{proof}
    It is easy to see $f$ maps the units of the operad $\termop_{\OwR}$ to the units of $\A(\opcat{O}) \BV \A(\opcat{R})$. To show that $f$ also preserves the composition, we first study an example that is small enough with respect to its data, yet not trivial.

    Let $X \in \opcat{O}$ with $Ar_{\opcat{O}}(X) = (c_1, c_2, c_3, c_4; c)$, for some $c_1, c_2, c_3, c_4, c \in \pi_0\opcat{O}$, and $Y_1, Y_2, Y_3, Y_4 \in \opcat{R}$ with $Ar(Y_{1 \leq i \leq 4}) = (d_1^i, d_2^i; d)$, for some $d_1^i, d_2^i, d \in \pi_0\opcat{R}$. Then $(X; Y_1, Y_2, Y_3, Y_4)$ is an object of the category $\OwR$. Similarly, let $z \in \opcat{O}$ with $Ar(z) = (t_1, t_2; c)$, for some $t_1,t_2 \in \pi_0\opcat{O}$, and $W_1, W_2 \in \opcat{R}$ with $Ar(w_i)=(r_1^i, r_2^i;d)$, for some $r_1^i, r_2^i \in \pi_0\opcat{R}$. Then $(z; W_1, W_2)$ is also an object in $\OwR$.

    Suppose there is a morphism
    \[(\varphi, \Phi): (x; Y_1, Y_2, Y_3, Y_4) \to (z; W_1, W_2)\] in $\OwR$, given by $\varphi: x \to z$ in $\opcat{O}$, with arity
    \[\begin{tikzcd}[ampersand replacement=\&,cramped]
    	{c_1} \& {c_2} \& {c_3} \& {c_4} \\
    	\& {t_1} \& {t_2}
    	\arrow[maps to, from=1-1, to=2-3]
    	\arrow[maps to, from=1-2, to=2-3]
    	\arrow[maps to, from=1-3, to=2-2]
    	\arrow[maps to, from=1-4, to=2-2]
    \end{tikzcd}\]
and fibers $\varphi^{-1}(1) = E_1$, $\varphi^{-1}(2) = E_2$. We compute $Ar(E_1)$ in detail. The computation of the arities of the other fiber objects in this example will be similar and therefore the details will be omitted.

As in Observation \ref{obs-fib}, $E_1$ belongs to the same connected component as $id_z^{-1}(1)$, hence $Ar(E_1)$ has $t_1$ as its root color. To find the color corresponding to $id_{E_1}^{-1}(1)$, we apply axiom (AX4) to the following commuting triangle
\[
\begin{tikzcd}[ampersand replacement=\&,cramped, column sep=tiny]
	{id_{E_1}^{-1}(1)} \& {E_1} \&\& {E_1} \\
	{id_X^{-1}(3)} \& x \&\& x \\
	\&\& z
	\arrow["\rhd"{marking, allow upside down, pos=0.3}, draw=none, from=1-1, to=1-2]
	\arrow["{=}"{marking}, draw=none, from=1-1, to=2-1]
	\arrow["{id_{E_1}}", from=1-2, to=1-4]
	\arrow["\rhd"{marking, allow upside down}, draw=none, from=1-2, to=2-2]
	\arrow["\rhd"{marking, allow upside down}, draw=none, from=1-4, to=2-4]
	\arrow["\rhd"{marking, allow upside down, pos=0.2}, draw=none, from=2-1, to=2-2]
	\arrow["id_X", from=2-2, to=2-4]
	\arrow["\varphi"', from=2-2, to=3-3]
	\arrow["\varphi", from=2-4, to=3-3]
\end{tikzcd}.
\]
Therefore $Ar(E_1) = (c_3, c_4; t_1)$ and similarly, $Ar(E_2) = (c_1, c_2; t_2)$

The cardinality $\card{\varphi}$ determines the domains and codomains of the morphisms in the family $\Phi = \{\varphi_{12}, \varphi_{22}, \varphi_{31}, \varphi_{41}\}$. Suppose the arities of the moprhisms $\varphi_{ij}$ are
\[\arraycolsep=10pt\def\arraystretch{1.6} 
\begin{array}{cccc}
    Ar(\varphi_{12}): & Ar(\varphi_{22}): & Ar(\varphi_{31}):& Ar(\varphi_{41}):\\
\begin{tikzcd}[ampersand replacement=\&,cramped]
	{d_1^1} \& {d^1_2} \\
	{r^2_1} \& {r^2_2}
	\arrow[maps to, from=1-1, to=2-2]
	\arrow[maps to, from=1-2, to=2-1]
\end{tikzcd} 
& 
\begin{tikzcd}[ampersand replacement=\&,cramped]
	{d^2_1} \& {d^2_2} \\
	{r^2_1} \& {r^2_2}
	\arrow[maps to, from=1-1, to=2-2]
	\arrow[maps to, from=1-2, to=2-1]
\end{tikzcd} 
& 
\begin{tikzcd}[ampersand replacement=\&,cramped]
	{d^3_1} \& {d^3_2} \\
	{r^1_1} \& {r^1_2}
	\arrow[maps to, from=1-1, to=2-2]
	\arrow[maps to, from=1-2, to=2-1]
\end{tikzcd}  
& 
\begin{tikzcd}[ampersand replacement=\&,cramped]
	{d^4_1} \& {d^4_2} \\
	{r^1_1} \& {r^1_2}
	\arrow[maps to, from=1-1, to=2-2]
	\arrow[maps to, from=1-2, to=2-1]
\end{tikzcd}  
\end{array}
\]

Denote the fibers of $\varphi_{ij}$
\[\arraycolsep=10pt\def\arraystretch{1.6} 
\begin{array}{cccc}
    \varphi_{12}^{-1}(1) = S^{12}_1 & \varphi_{22}^{-1}(1) = S^{22}_1 & \varphi_{31}^{-1}(1) = S^{31}_1 & \varphi_{41}^{-1}(1) = S^{41}_1 \\
    \varphi_{12}^{-1}(2) = S^{12}_2 & \varphi_{22}^{-1}(2) = S^{22}_2 & \varphi_{31}^{-1}(2) = S^{31}_2 & \varphi_{41}^{-1}(2) = S^{41}_2 
\end{array}.
\] 
We compute $Ar(S^{ij}_1) = (d^i_2; r^{j}_1)$ and $Ar(S^{ij}_2) = (d^i_1; r^{j}_2)$.

The arity of $(\varphi, \Phi)$ is 
\[\begin{tikzcd}[ampersand replacement=\&,cramped,column sep=small]
	{(c_1, d^1_1)} \& {(c_1, d^1_2)} \& {(c_2, d^2_1)} \& {(c_2, d^2_2)} \& {(c_3, d^3_1)} \& {(c_3, d^3_2)} \& {(c_4, d^4_1)} \& {(c_4, d^4_2)} \\
	\\
	\\
	\&\& {(t_1,r_1^1)} \& {(t_1, r^1_2)} \& {(t_2, r^2_1)} \& {(t_2, r^2_2)}
	\arrow[draw={rgb,255:red,153;green,92;blue,214}, maps to, from=1-1, to=4-6, line width=0.2mm]
	\arrow[draw={rgb,255:red,153;green,92;blue,214}, maps to, from=1-2, to=4-5, line width=0.2mm]
	\arrow[draw={rgb,255:red,92;green,214;blue,92}, curve={height=-6pt}, maps to, from=1-3, to=4-6, line width=0.2mm]
	\arrow[draw={rgb,255:red,92;green,214;blue,92}, maps to, from=1-4, to=4-5, line width=0.2mm]
	\arrow[draw={rgb,255:red,51;green,146;blue,255}, curve={height=-12pt}, maps to, from=1-5, to=4-4, line width=0.2mm]
	\arrow[draw={rgb,255:red,51;green,146;blue,255}, maps to, from=1-6, to=4-3, line width=0.2mm]
	\arrow[draw={rgb,255:red,214;green,92;blue,92}, maps to, from=1-7, to=4-4, line width=0.2mm]
	\arrow[draw={rgb,255:red,214;green,92;blue,92}, maps to, from=1-8, to=4-3, line width=0.2mm]
\end{tikzcd}.\]
and its four fibers of $(\varphi, \Phi)$ are:
    \[\arraycolsep=10pt\def\arraystretch{1.6}
    \begin{array}{cc}
         (\varphi, \Phi)^{-1}(1)=(E_1; S^{31}_{1}, S^{41}_{1}), & (\varphi, \Phi)^{-1}(2)=(E_1; S^{31}_{2}, S^{41}_{2}), \\
         (\varphi, \Phi)^{-1}(3) = (E_2; S^{12}_{1}, S^{22}_{1}), & 
         (\varphi, \Phi)^{-1}(4) = (E_2; S^{12}_{2}, S^{22}_{2})
    \end{array}
    \]

    The data above implies that equalities hold:
    \begin{itemize}  \setlength\itemsep{1em}
        \item $\mu_{\termop_{\OwR}}^{(\varphi, \Phi)}
        \Big(
        (\varphi, \Phi)^{-1}(i)_{1 \leq i \leq 4},
        (Z; W_1, W_2)
        \Big)
        = (X; Y_1, Y_2, Y_3, Y_4)$ in $\termop_{\OwR}$;
        \item $\mu_{\A(\opcat{O})}^{Ar \varphi} ([E_1], [E_2], [Z]) = [X]$ in $\A(\opcat{O})$, and therefore, \\
        $\mu_{BV}^{Ar \varphi} ([E_1] \otimes d, [E_2] \otimes d, [Z] \otimes d) = [X] \otimes d$ in $\A(\opcat{O}) \BV \A(\opcat{R})$, for any color $d \in \pi_0\opcat{R}$;
        \item $\mu_{\A(\opcat{R})}^{Ar \varphi_{ij}}([S^{ij}_1], [S^{ij}_2], [w_j]) = [y_i]$ in $\A(\opcat{R})$, and therefore, \\
        $\mu_{BV}^{Ar \varphi_{ij}}(c \otimes [S^{ij}_1], c \otimes [S^{ij}_2], c \otimes [w_j]) = c \otimes [y_i]$, for any color $c \in \pi_0\opcat{O}$.
    \end{itemize}

We wish to show that 
\begin{equation} \label{exp-first}
\arraycolsep=0pt\def\arraystretch{1.6}
\begin{array}{ccll}
     \mu_{BV}^{Ar(\varphi, \Phi)} & \Big( &  \mu_{BV}^{\nu}(c_3 \otimes [S^{31}_1], c_4 \otimes [S^{41}_1], [E_1] \otimes r^1_1), &\\
     & &  \mu_{BV}^{\nu}(c_3 \otimes [S^{31}_2], c_4 \otimes [S^{41}_2],  [E_1] \otimes r^1_2), &\\
     & &  \mu_{BV}^{\nu}(c_1 \otimes [S^{12}_1], c_2 \otimes [S^{22}_1], [E_2] \otimes r^2_1), &\\
     & &  \mu_{BV}^{\nu}(c_1 \otimes [S^{12}_2], c_2 \otimes [S^{22}_2],  [E_2] \otimes r^2_2), &\\
     & & \mu_{BV}^{\nu}(t_1 \otimes [W_1], t_2 \otimes [W_2], [Z] \otimes d) & \Big)
\end{array} 
\end{equation}
is equal to 
\[\mu_{BV}^{\nu}(c_1 \otimes [Y_1], c_2 \otimes [Y_2], c_3 \otimes [Y_3], c_4 \otimes [Y_4], [X] \otimes d).\]
We first decompose the morphism $Ar(\varphi, \Phi)$ into the composite
\[\begin{tikzcd}[ampersand replacement=\&,cramped,column sep=tiny]
	\& {(c_1, d^1_1)} \& {(c_1, d^1_2)} \& {(c_2, d^2_1)} \& {(c_2, d^2_2)} \& {(c_3, d^3_1)} \& {(c_3, d^2_2)} \& {(c_4, d^4_1)} \& {(c_4, d^4_2)} \\
	{\coprod (c_i, Ar(\varphi_{ij}))} \\
	\& {(c_1, r^2_1)} \& {(c_1, r^2_2)} \& {(c_2, r^2_1)} \& {(c_2, r^2_2)} \& {(c_3, r^1_1)} \& {(c_3, r^1_2)} \& {(c_4, r^1_1)} \& {(c_4, r^1_2)} \\
	{\rho:} \\
	\\
	\&\& {(t_1, r^1_1)} \& {(t_1, r^1_2)} \&\&\& {(t_2,r^2_1)} \& {(t_2, r^2_2)}
	\arrow[draw={rgb,255:red,153;green,92;blue,214}, maps to, from=1-2, to=3-3, line width=0.2mm]
	\arrow[draw={rgb,255:red,153;green,92;blue,214}, maps to, from=1-3, to=3-2, line width=0.2mm]
	\arrow[draw={rgb,255:red,92;green,214;blue,92}, maps to, from=1-4, to=3-5, line width=0.2mm]
	\arrow[draw={rgb,255:red,92;green,214;blue,92}, maps to, from=1-5, to=3-4, line width=0.2mm]
	\arrow[draw={rgb,255:red,51;green,102;blue,255}, maps to, from=1-6, to=3-7, line width=0.2mm]
	\arrow[draw={rgb,255:red,51;green,102;blue,255}, maps to, from=1-7, to=3-6, line width=0.2mm]
	\arrow[draw={rgb,255:red,214;green,92;blue,92}, maps to, from=1-8, to=3-9, line width=0.2mm]
	\arrow[draw={rgb,255:red,214;green,92;blue,92}, maps to, from=1-9, to=3-8, line width=0.2mm]
	\arrow[draw={rgb,255:red,153;green,92;blue,214}, curve={height=6pt}, maps to, from=3-2, to=6-7, line width=0.2mm]
	\arrow[draw={rgb,255:red,153;green,92;blue,214}, curve={height=6pt}, maps to, from=3-3, to=6-8, line width=0.2mm]
	\arrow[draw={rgb,255:red,92;green,214;blue,92}, maps to, from=3-4, curve={height=-6pt}, to=6-7, line width=0.2mm]
	\arrow[draw={rgb,255:red,92;green,214;blue,92}, maps to,curve={height=-6pt},  from=3-5, to=6-8, line width=0.2mm]
	\arrow[draw={rgb,255:red,51;green,102;blue,255}, maps to, curve={height=6pt}, from=3-6, to=6-3, line width=0.2mm]
	\arrow[draw={rgb,255:red,51;green,102;blue,255}, maps to, curve={height=6pt}, from=3-7, to=6-4, line width=0.2mm]
	\arrow[draw={rgb,255:red,214;green,92;blue,92}, curve={height=-6pt}, maps to, from=3-8, to=6-3, line width=0.2mm]
	\arrow[draw={rgb,255:red,214;green,92;blue,92}, curve={height=-6pt}, maps to, from=3-9, to=6-4, line width=0.2mm]
\end{tikzcd}\]
Since the axiom (COMP) holds, the expression (\ref{exp-first}) is equal to
\begin{equation} \label{exp-comp}
\arraycolsep=0pt\def\arraystretch{1.6}
\begin{array}{ccllrll} 
     \mu^{\coprod (c_i, Ar(\varphi_{ij}))}_{BV} & \Big( &
            c_1 \otimes [S^{12}_1], c_1 \otimes [S^{12}_2], & & &\\
     & &    c_2 \otimes [S^{22}_1], c_2 \otimes [S^{22}_2], & & &\\
     & &    c_3 \otimes [S^{31}_1], c_3 \otimes [S^{31}_2], & & &\\
     & &    c_4 \otimes [S^{41}_1], c_4 \otimes [S^{41}_2],  &\mu_{BV}^{\rho}\big( & [E_1] \otimes r^1_1, [E_1] \otimes r^1_2,
                             [E_2] \otimes r^2_1, [E_2] \otimes r^2_2, \\
    &&&&\mu_{BV}^{\nu}(t_1 \otimes [W_1], t_2 \otimes [W_2], [Z] \otimes d) &\big) & \Big)
\end{array} 
\end{equation} 
Similarly, due to (COMP), the following equality holds
\begin{equation}\label{exp-rho1-rho2}
    \arraycolsep=0pt\def\arraystretch{1.6}
    \begin{array}{crc}
         \mu_{BV}^{\rho}\big(& [E_1] \otimes r^1_1,[E_1] \otimes r^1_2,
                             [E_2] \otimes r^2_1, [E_2] \otimes r^2_2,\\
         &                  \mu_{BV}^{\nu}(t_1 \otimes [W_1], t_2 \otimes [W_2], [Z] \otimes d) \big)&= \\
         =\mu_{BV}^{\nu \circ \rho}  \big(& \mu^{\rho_1}_{BV}([E_1] \otimes r^1_1, [E_1] \otimes r^1_2, t_1 \otimes [W_1]), \\
         &                                  \mu^{\rho_2}_{BV}([E_2] \otimes r^2_1, [E_2] \otimes r^2_2, t_2 \otimes [W_2]), \\
         &[Z] \otimes d\big).
    \end{array}
\end{equation}
The morphisms $\rho_1$ and $\rho_2$ are given by
\[\begin{tikzcd}[ampersand replacement=\&,cramped,column sep=tiny]
	\& {(c_3, r^1_1)} \& {(c_3, r^1_2)} \& {(c_4, r^1_1)} \& {(c_4, r^1_2)} \& {(c_1, r^2_1)} \& {(c_1, r^2_2)} \& {(c_2, r^2_1)} \& {(c_2, r^2_2)} \& \\
	{\rho_1:} \&\&\&\&\&\&\&\&\& {:\rho_2.} \\
	\\
	\&\& {(t_1, r^1_1)} \& {(t_1, r^1_2)} \&\&\& {(t_2,r^2_1)} \& {(t_2, r^2_2)}
	\arrow[draw={rgb,255:red,51;green,102;blue,255}, curve={height=6pt}, maps to, from=1-2, to=4-3, line width=0.2mm]
	\arrow[draw={rgb,255:red,51;green,102;blue,255}, curve={height=6pt}, maps to, from=1-3, to=4-4, line width=0.2mm]
	\arrow[draw={rgb,255:red,214;green,92;blue,92}, curve={height=-6pt}, maps to, from=1-4, to=4-3, line width=0.2mm]
	\arrow[draw={rgb,255:red,214;green,92;blue,92}, curve={height=-6pt}, maps to, from=1-5, to=4-4, line width=0.2mm]
	\arrow[draw={rgb,255:red,153;green,92;blue,214}, curve={height=6pt}, maps to, from=1-6, to=4-7, line width=0.2mm]
	\arrow[draw={rgb,255:red,153;green,92;blue,214}, curve={height=6pt}, maps to, from=1-7, to=4-8, line width=0.2mm]
	\arrow[draw={rgb,255:red,92;green,214;blue,92}, curve={height=-6pt}, maps to, from=1-8, to=4-7, line width=0.2mm]
	\arrow[draw={rgb,255:red,92;green,214;blue,92}, curve={height=-6pt}, maps to, from=1-9, to=4-8, line width=0.2mm]
\end{tikzcd}\]
We now focus on the expression $\mu^{\rho_1}_{BV}([E_1] \otimes r^1_1, [E_1] \otimes r^1_2, t_1 \otimes [W_1])$. We rewrite it in the language of classical colored operads
\[\mu^{\rho_1}_{BV}([E_1] \otimes r^1_1, [E_1] \otimes r^1_2, t_1 \otimes [W_1]) =  \gamma_{BV}(t_1 \otimes [W_1], [E_1] \otimes r^1_1, [E_1] \otimes r^1_2) \cdot \pi(\rho_1),\]
where $\pi(\rho_1)$ is the permutation in the decomposition of $\rho_1$. Moreover, $\pi(\rho_1)$ is the permutation that reorders colors from lexicographical notation to reverse lexicographical notation, i.e., \textit{shuffle} in the sense of Definition \ref{def-bv}. Hence, 
\begin{equation} \label{exp-interchnage}
\arraycolsep=0pt\def\arraystretch{1.6}
\begin{array}{ccc}
    \mu^{\rho_1}_{BV}([E_1] \otimes r^1_1, [E_1] \otimes r^1_2, t_1 \otimes [W_1]) &=& \gamma_{BV}([E_1] \otimes d, c_3 \otimes [W_1], c_4 \otimes [W_1])\\
    &=& \mu^{\nu}_{BV}(c_3 \otimes [W_1], c_4 \otimes [W_1], [E_1] \otimes d).
\end{array}
\end{equation}
We apply the above to the righthand side of (\ref{exp-rho1-rho2}) to obtain equality 
\begin{equation}\label{exp-after-inter}
\arraycolsep=0pt\def\arraystretch{1.6}
\begin{array}{crc}
      \mu_{BV}^{\rho}\big( & [E_1] \otimes r^1_1, [E_1] \otimes r^1_2,
                             [E_2] \otimes r^2_1, [E_2] \otimes r^2_2, \\
                           & \mu_{BV}^{\nu}(t_1 \otimes [W_1], t_2 \otimes [W_2], [Z] \otimes d) \big) & = \\
    = \mu_{BV}^{\nu \circ \rho}  \big( &     \mu^{\nu}_{BV}(c_3 \otimes [W_1], c_4 \otimes [W_1], [E_1] \otimes d), \\
                                        & \mu^{\nu}_{BV}(c_1 \otimes [W_2], c_2 \otimes [W_2], [E_2] \otimes d), \\
                                        & [Z] \otimes d\big).
\end{array}
\end{equation}
The composite $\nu \circ \rho$ 
is equal to the composite $(Ar\varphi \otimes d) \circ \nu$, 
\[\begin{tikzcd}[ampersand replacement=\&,cramped,cramped,column sep=tiny]
	\& {(c_1, r^2_1)} \& {(c_1, r^2_2)} \& {(c_2, r^2_1)} \& {(c_2, r^2_2)} \& {(c_3, r^1_1)} \& {(c_3, r^1_2)} \& {(c_4, r^1_1)} \& {(c_4, r^1_2)} \\
	{\nu:} \\
	\\
	\&\& {(c_1, d)} \& {(c_2, d)} \&\&\& {(c_3,d)} \& {(c_4, d)} \\
	{Ar \varphi \otimes d:} \\
	\&\&\& {(t_1, d)} \&\&\& {(t_2, d)}
	\arrow[draw={rgb,255:red,153;green,92;blue,214}, curve={height=6pt}, maps to, from=1-2, to=4-3, line width=0.2mm]
	\arrow[draw={rgb,255:red,153;green,92;blue,214}, curve={height=6pt}, maps to, from=1-3, to=4-3, line width=0.2mm]
	\arrow[draw={rgb,255:red,92;green,214;blue,92}, curve={height=6pt}, maps to, from=1-4, to=4-4, line width=0.2mm]
	\arrow[draw={rgb,255:red,92;green,214;blue,92}, curve={height=6pt}, maps to, from=1-5, to=4-4, line width=0.2mm]
	\arrow[draw={rgb,255:red,51;green,102;blue,255}, curve={height=-6pt}, maps to, from=1-6, to=4-7, line width=0.2mm]
	\arrow[draw={rgb,255:red,51;green,102;blue,255}, curve={height=-6pt}, maps to, from=1-7, to=4-7, line width=0.2mm]
	\arrow[draw={rgb,255:red,214;green,92;blue,92}, curve={height=-6pt}, maps to, from=1-8, to=4-8, line width=0.2mm]
	\arrow[draw={rgb,255:red,214;green,92;blue,92}, curve={height=-6pt}, maps to, from=1-9, to=4-8, line width=0.2mm]
	\arrow[curve={height=6pt}, shift right=3, maps to, from=4-3, to=6-7, line width=0.2mm]
	\arrow[curve={height=6pt}, maps to, from=4-4, to=6-7, line width=0.2mm]
	\arrow[curve={height=-6pt}, maps to, from=4-7, to=6-4, line width=0.2mm]
	\arrow[curve={height=-6pt},  shift left=3, maps to, from=4-8, to=6-4, line width=0.2mm]
\end{tikzcd}.\]
We apply this to the right-hand side of (\ref{exp-after-inter}) together with  axiom (COMP),
\begin{equation}
\arraycolsep=0pt\def\arraystretch{1.6}
    \begin{array}{lrccc} 
    \mu_{BV}^{\nu \circ \rho}  \big(& \mu^{\nu}_{BV}(c_3 \otimes [W_1], c_4 \otimes [W_1], [E_1] \otimes d), \\
            &\mu^{\nu}_{BV}(c_1 \otimes [W_2], c_2 \otimes [W_2], [E_2] \otimes d),\\
            &[Z]\otimes d\big)&=\\
         =\mu_{BV}^{(Ar\varphi \otimes d) \circ \nu}  \big( & \mu^{\nu}_{BV}(c_3 \otimes [W_1], c_4 \otimes [W_1], [E_1] \otimes d),\\
        &\mu^{\nu}_{BV}(c_1 \otimes [W_2], c_2 \otimes [W_2], [E_2] \otimes d),\\
        & [Z]\otimes d\big) &=\\
        = \mu^{\nu}_{BV}(&c_1 \otimes [W_2], c_2 \otimes [W_2], c_3 \otimes [W_1], c_4 \otimes [W_1], \\
        &\mu^{Ar \varphi \otimes d}_{BV} ([E_1] \otimes d, [E_2] \otimes d, [Z] \otimes d)) &= \\
        = \mu^{\nu}_{BV}(&c_1 \otimes [W_2], c_2 \otimes [W_2], c_3 \otimes [W_1], c_4 \otimes [W_1], \\
        &[X] \otimes d). 
    \end{array}
\end{equation}
We insert this result back to (\ref{exp-comp}) and apply (COMP)
\begin{equation}\label{exp-final}
\arraycolsep=0pt\def\arraystretch{1.6}
\begin{array}{lcllrll} 
     \mu^{\coprod (c_i \otimes Ar(\varphi_{ij}))}_{BV} & \Big( &
            c_1 \otimes [S^{12}_1], c_1 \otimes [S^{12}_2], c_2 \otimes [S^{22}_1], c_2 \otimes [S^{22}_2], & & &\\
     & &    c_3 \otimes [S^{31}_1], c_3 \otimes [S^{31}_2], c_4 \otimes [S^{41}_1], c_4 \otimes [S^{41}_2],  \\
     & & \mu^{\nu}_{BV}(c_1 \otimes [W_2], c_2 \otimes [W_2], c_3 \otimes [W_1], c_4 \otimes [W_1], [X] \otimes d) \Big)= \\
     =\mu^{\nu \circ \coprod (c_i \otimes Ar(\varphi_{ij}))}_{BV} & \Big( &
        \mu_{BV}^{c_1 \otimes Ar\varphi_{12}}(c_1 \otimes [S^{12}_1], c_1 \otimes [S^{12}_2], c_1 \otimes [W_2]), \\ 
    &&  \mu_{BV}^{c_2 \otimes Ar\varphi_{22}}(c_2 \otimes [S^{22}_1], c_2 \otimes [S^{22}_2], c_2 \otimes [W_2]), \\
    &&  \mu_{BV}^{c_3 \otimes Ar\varphi_{31}}(c_3 \otimes [S^{31}_1], c_3 \otimes [S^{31}_2], c_3 \otimes [W_1]), \\
    &&  \mu_{BV}^{c_4 \otimes Ar\varphi_{41}}(c_4 \otimes [S^{41}_1], c_4 \otimes [S^{41}_2], c_4 \otimes [W_1]), [X] \otimes d \Big) \\
    =\mu^{\nu \circ \coprod (c_i \otimes Ar(\varphi_{ij}))}_{BV} & \Big( &
        c_1 \otimes [Y_1], c_2 \otimes [Y_2], c_3 \otimes [Y_3], c_4 \otimes [Y_4], [X] \otimes d \Big)    
\end{array} 
\end{equation} 
Finally, the composing $(c_i \otimes Ar(\varphi_{ij}))$ with an order preserving $\nu$ ignores the permutations performed inside the direct summand. Therefore, the latter expression is just \[\mu_{BV}^{\nu}(c_1 \otimes [Y_1], c_2 \otimes [Y_2], c_3 \otimes [Y_3], c_4 \otimes [Y_4], [X] \otimes d).\]
The procedure above can be used for any composition in $\termop_{\OwR}$, the only complication is a more detailed bookkeeping of indexes. We conclude that 
\[
\begin{array}{rrcl}
    f_{\opcat{O}, \opcat{R}}: & \termop_{\OwR} &  \to & Ar_{\OwR}^{\ast}(\A(\opcat{O}) \BV \A(\opcat{R}))  \\
     &  (x; Y_1, \ldots, Y_n) &\mapsto& \mu_{BV}^{\nu}\left(\left(c_1 \otimes [Y_1], \ldots, c_n \otimes [Y_n]\right), [X] \otimes d\right)
\end{array}
\]
is indeed a morphism of operads.
\end{proof}

The morphisms $f_{\opcat{O}, \opcat{R}}$ compose into a natural transformation $f: \termop_{-\wr-} \Rightarrow \A(-)\BV\A(-)$. It is a basic exercise in category theory that the corresponding morphisms $\gamma_{\opcat{O}, \opcat{R}}$ compose into a natural transformation 
$\gamma: \A(- \wr -) \Rightarrow \A(-) \BV \A(-)$.

\begin{prop}
    Given operadic categories $\opcat{O}, \opcat{R}$, the morphism $\gamma_{\opcat{O}, \opcat{R}}: \A(\OwR) \to \A(\opcat{O}) \BV \A(\opcat{R})$ is an isomorphism. It's inverse is given by 
    \begin{equation*}
        \begin{array}{crcl}
             \psi_{\opcat{O}, \opcat{R}}:& \A(\opcat{O}) \BV \A(\opcat{R}) &\to & \A(\OwR)  \\
             & [X] \otimes d & \mapsto & [(X; V_d, \ldots, V_d)] \\
             & c \otimes [Y] & \mapsto & [(U_c; Y)]
        \end{array};
    \end{equation*}
    where $X \in \opcat{O}$, $Y \in \opcat{R}$, $U_c$ (resp. $V_d$) is the chosen local terminal objects corresponding to $c \in \pi_0\opcat{O}$ (resp. $d \in \pi_0\opcat{R}$).
\end{prop}
\begin{proof}
We check that $\psi_{\opcat{O}, \opcat{R}}$ is a well-defined morphism of operads. Fix a color $d \in \pi_0\opcat{R}$ and a corresponding local terminal object $V_d \in \opcat{R}$. There is a functorial inclusion $i_{d}: \opcat{O} \to \OwR$ that acts on objects by $i_d(x) = (x; V_d, \ldots, V_d)$. We recall the unit $\eta: Id_{\Op} \Rightarrow inc \circ \Arleft$ of the adjunction  $\xymatrix{\Arleft:\Op\ar[r]<1ex> & \SOp:inc\ar[l]<1ex>_(0.55){\upvdash}}$. The square
\[\begin{tikzcd}[ampersand replacement=\&,cramped]
	{\termop_{\opcat{O}}} \& {\termop_{\OwR}} \\
	{\A(\opcat{O})} \& {\A(\OwR)}
	\arrow["{\termop_{i_d}}", from=1-1, to=1-2]
	\arrow["{\eta_{\termop_{\opcat{O}}}}"', from=1-1, to=2-1]
	\arrow["{\eta_{\termop_{\OwR}}}", from=1-2, to=2-2]
	\arrow["{\A(i_d)}"', from=2-1, to=2-2]
\end{tikzcd}\]
commutes in $\Op$ and therefore $\psi_{\opcat{O}, \opcat{R}}$ preserves the composition of generators of the type $[X]\otimes d$, for any fixed $d$ and $[X] \in \A(\opcat{O})$. To show that $\psi_{\opcat{O}, \opcat{R}}$ preserves the composition of $c \otimes \A(\opcat{R})$ for any $c \in \pi_0\opcat{O}$, it is enough to apply the same argument to the functorial inclusion $j_c: \opcat{R} \to \OwR$ given on objects by $(U_c; Y)$, where $U_c$ is the chosen local terminal object corresponding to $c$.

To see that $\psi_{\opcat{O}, \opcat{R}}$ preserves the interchange relation, fix $x \in \opcat{O}$ and $y \in \opcat{R}$. Denote $c := \pi_0(x)$, $d := \pi_0(Y)$, and $U_c, V_d$ as the corresponding local terminal objects in their respective categories. There is a morphism 
\[(id_X, !_Y): (X; Y, \ldots, Y) \to (X; V_d, \ldots, V_d)\] in $\OwR$, where $!_Y$ is a family of terminal morphisms $!_Y: Y \to V_d$. For any $i \in \card{(X; V_d, \ldots, V_d)} = \card{X}$, the $i$-th fiber $(id_X, !_Y)^{-1}(i)$ is $(U_{i}; Y)$, where $U_i$ is the chosen local terminal object of $\pi_0(id_X^{-1}(i))$.

Similarly, there is a morphism
\[(!_X, id_Y): (X; Y, \ldots, Y) \to (U_c; Y)\] in $\OwR$, where $id_Y$ is a family of the identity morphisms. For any $i \in \card{(U_c; Y)} = \card{y}$, the $i$-th fiber $(!_X, id_Y)^{-1}(i)$ is $(X; V_j, \ldots, V_j)$, where $V_j$ is the chosen local terminal object of $\pi_0(id_Y^{-1}(j))$.

By construction of $\A$, the following equality holds
\begin{multline}
\arraycolsep=0pt\def\arraystretch{1.6}
    \mu_{\A(\OwR)}^{Ar(id_X, !_Y)}([(U_{1}; Y)], \ldots, [(U_{\card{X}}; Y)], [(X; V_d, \ldots, V_d)]) = [(X; Y, \ldots, Y)] =\\
         = \mu_{\A(\OwR)}^{Ar(!_X, id_Y)}([(X; V_1, \ldots, V_1)], \ldots, [(X; V_{\card{y}}, \ldots, V_{\card{y}})], [(U_c; Y)]).
\end{multline}
holds in the operad $\A(\OwR)$. We observe that $Ar(id_X, !_Y)$ is the composition $\textit{shuffle} \circ Ar(!_X, id_Y)$ and conclude that $\psi_{\opcat{O}, \opcat{R}}$ preserves the interchange relation. Therefore, $\psi_{\opcat{O}, \opcat{R}}$ is a well-defined morphism of operads, and it is easy to see that it is the inverse of $\gamma_{\opcat{O}, \opcat{R}}$.
\end{proof}

We have so far defined the isomorphism $\iota: \A(\mathbf{1}) \xrightarrow{\sim} \star$ and the natural isomorphism $\gamma: \A(- \wr -) \Rightarrow \A(-) \BV \A(-)$. It is straightforward to show that $\iota$ and $\gamma$ satisfy the known associativity and unitality conditions. We therefore conclude the main result of the present section.
\begin{thm} \label{thm-mon-fun}
    The functor $\A: \CatOp \to \SOp$ is a strong monoidal functor.
\end{thm}
\begin{cor} \label{cor-main}
    For  a $\frC$-colored unital symmetric operad $\oper{P}$ and a $\frD$-colored unital symmetric operad $\oper{Q}$,
    \[\oper{P} \BV \oper{Q} \cong \A(\I\oper{P} \wr \I\oper{Q}).\]
\end{cor}
\begin{proof}
    Apply Theorem \ref{thm-mon-fun} and Observation \ref{obs-AI-op}.
\end{proof}

%% file: relation.tex
We compare the wreath product of operadic categories with two other types of wreath product: the categorical wreath product, as defined by C.~Berger, and the wreath product of operator categories, as defined by C.~Barwick.

\subsection*{Categorical wreath product}

C.~Berger defines, for each small category $\mathcal{A}$, the categorical wreath product $\Delta \wr_{\text{B}} \mathcal{A}$ and $\Gamma \wr_{\text{B}} \mathcal{A}$. The categorical wreath product is then applied to construct a category~$\Theta^n$; this category has a remarkable property that the homotopy theory of $n$-fold loop spaces is equivalent to the homotopy theory of reduced  $\Theta^n$-spaces\cite{BERGER2007230}. In this section, we recollect relevant notions from the original paper and explain the relation between the categorical wreath product and the operadic one. 

Recall that the category $\Delta$ is the category of finite non-empty ordinals $[n] = \{0 \leq \ldots \leq n\}$ together with order-preserving maps in the notation. The category $\Gamma$, introduced by Segal \cite[Def.~1.1.]{Segal}, is the category of finite sets $\underline{n} = \{1 \leq \ldots \leq n\}$ with morphisms given by ordered $m$-tuples of pairwise disjoint subsets of $n$. Thus defined, the category $\Gamma$ is isomorphic to the opposite category of finite pointed sets $\Fin_*^{op}$, morphisms in $\Gamma$ are described by preimages of morphisms in $\Fin_{\ast}$.  Segal also defines a functor $\gamma: \Delta \to \Gamma$ which takes the ordinal $[m]$ to $\underline{m}$ and an order-preserving morphism $\varphi: [m] \to [n]$ to a morphism $\gamma(\varphi): \underline{m} \to \underline{n}$ defined by $\gamma(\varphi)(i) = \{\varphi(i-1) + 1 \leq \ldots \leq \varphi(i)\}$, for $i \in \underline{m}$.

\begin{defn}[\cite{BERGER2007230} Def.~3.1.] \label{defn-berger-wr}
    The categorical wreath product $\Delta \wr_{\text{B}} \mathcal{A}$ (resp. $\Gamma \wr_{\text{B}} \mathcal{A}$) is defined to be the category
    \begin{itemize}
        \item with objects the $m$-tuples $(A_1, \ldots, A_m)$ of objects of $\mathcal{A}$, for varying $m \geq 0$;
        \item with morphisms $(\varphi; \varphi_1, \ldots, \varphi_m): (A_1, \ldots, A_m) \to (B_1, \ldots, B_n)$, where $\varphi: [m] \to [n]$ is a morphism in $\Delta$ and \[\varphi_i =\{\varphi_{ij}: A_i \to B_j \mid \varphi(i-1) + 1 \leq j \leq \varphi(i)\}\]
        \[(\text{resp. } \varphi_i =\{\varphi_{ij}: A_i \to B_j \mid  j \in \varphi(i)\})\]
    \end{itemize}
    Composition in $\Delta \wr_{\text{B}} \mathcal{A}$ (resp. $\Gamma \wr_{\text{B}} \mathcal{A}$) is then directly induced by composition in $\mathcal{A}$.
\end{defn}

Given a category $\mathcal{B}$ equipped with a functor $\mathcal{B} \to \Gamma$, Definition~\ref{defn-berger-wr} can be straightforwardly generalized to an arbitrary categorical wreath product $\mathcal{B} \wr \mathcal{A}$. The wreath product $\Delta \wr_{\text{B}} \mathcal{A}$ is then given by $\gamma: \Delta \to \Gamma$ and $\Gamma \wr_{\text{B}} \mathcal{A}$ is given by the identity functor \cite[Rem.~3.4.]{BERGER2007230} 

An operadic category $\opcat{O}$ is equipped with the cardinality functor $\card{-}: \opcat{O} \to \Fin$, which we compose with the inclusion $\Fin \hookrightarrow \Fin_{\ast}$. The opposite of this composition defines a functor $\opcat{O}^{op} \to \Gamma$.

\begin{prop} \label{prop-bg-bt}
    Suppose $\opcat{O}$ and $\opcat{R}$ are operadic categories, and $\card{-}^{op}: \opcat{O}^{op} \to \Gamma$ is as above. Then the categorical wreath product $\opcat{O}^{op} \wr_{\text{B}} \opcat{R}^{op}$ in the sense of Definition \ref{defn-berger-wr} is isomorphic to the opposite category of the operadic wreath product $\OwR^{op}$ in the sense of Definition \ref{defn-prewr}. 
\end{prop}
\begin{proof}
    The objects of the categorical wreath product $\opcat{O}^{op} \wr_{\text{B}} \opcat{R}^{op}$ are tuples $(A; B_1, \ldots, B_n)$, where $A \in \opcat{O}^{op}$, $\card{A}^{op} = \underline{n}$ and $B_1, \ldots, B_n \in \opcat{R}^{op}$. It is easy to see that the objects of $\opcat{O}^{op} \wr_{\text{B}} \opcat{R}^{op}$ coincide with the objects of the opposite category of the wreath product of operadic categories $\OwR^{op}$.

    A morphism $(\varphi; \varphi_1, \ldots, \varphi_m): (Z; W_1, \ldots w_k) \to (A; B_1,\ldots, B_n)$ in $\opcat{O}^{op} \wr \opcat{R}^{op}$ is given by a morphism $\varphi: Z \to A$ in $\opcat{O}^{op}$ (i.e., $\varphi: A \to Z$ in $\opcat{O}$) and a family of morphisms in $\opcat{R}^{op}$, for each $i \in \underline{k}$,
    \[\varphi_i = \{\varphi_{ij}: W_i \to B_j \mid j \in \card{\varphi}^{op}(i)\}.\]
    By construction, $\cup_{i \in \underline{k}} \card{\varphi}^{op}(i) = \underline{n}$, since the domain of $\card{\varphi} = \underline{n}$. Therefore, 
    \[\cup_{i \in \underline{k}} \varphi_i = \{\varphi_{ji} : B_j \to W_i \mid j \in \underline{n}, i = \card{\varphi}(j)\}.\] This is precisely the data of $\OwR$.
\end{proof}

Segal's functor $\gamma: \Delta \to \Gamma$, induces a functor $\gamma \wr \mathcal{A} \to \gamma \wr \Gamma$ which is the identity on objects, and sends the morphism $(\varphi; \varphi_1, \ldots, \varphi_m)$ in $\Delta \wr_{\text{B}} \mathcal{A}$ to a morphism $(\gamma(\varphi); \varphi_1, \ldots, \varphi_m)$ in  $\Gamma \wr_{\text{B}} \mathcal{A}$. Then there is a canonical assembly functor $\alpha: \Gamma \wr_{\text{B}} \Gamma \to \Gamma$ which takes $(\underline{n}_1, \ldots, \underline{n}_k)$ to $\underline{n_1 + \cdots + n_k}$ \cite[Lemma~3.2.]{BERGER2007230}.

Berger defines the category $\Theta_n$ together with canonical functors $\gamma_n: \Theta_n \to \Gamma$ as an iterated categorical wreath product of $\Delta$ with itself. Precisely $\Theta_1 := \Delta$ and $\gamma_1 := \gamma$; then $\Theta_{n} := \Delta \wr_{\text{B}} \Theta_{n-1}$ and $\gamma_n: \Theta_n \to \Gamma$ is the composition $\Delta \wr_{\text{B}} \Theta_{n-1} \xrightarrow{\gamma \wr \gamma_{n-1}} \Gamma \wr_{\text{B}} \Gamma \xrightarrow{\alpha} \Gamma$\cite[Def.~3.3.]{BERGER2007230}. 

Denote by $\Delta^{c}$ the subcategory of $\Delta$ generated by inner faces and degeneracies (such morphisms are called \textit{covers}) and consider the restriction $\gamma: \Delta^{c} \to \Gamma$.
\begin{obs}
    The restriction $\gamma: \Delta^{c} \to \Gamma$ factors through the inclusion $\Dalg^{op} \hookrightarrow \Gamma$. Moreover, it is an isomorphism $\Delta^{c} \cong \Dalg^{op}$.
\end{obs}
\begin{proof}
    Recall that objects in $\Dalg$ are denoted $\underline{n} = \{1, \ldots, n\}$ with $\underline{0} = \emptyset$. The the morphisms in $\Dalg^{op}$ are described by the preimages of order-preserving morphisms in $\Dalg$. By construction, the image of $\gamma: \Delta^{c} \to \Gamma$ lies in $\Dalg^{op}$ and is obviously bijective on objects. It is also easy to see that $\gamma$ is faithful. To see that it is full, let $f: \underline{n} \to \underline{m}$ be a morphism in $\Dalg^{op}$. We define $\hat{f}: [n] \to [m]$ iteratively.

    \begin{enumerate}
        \item Put $\hat{f}(0) := 0$ and let $last := f(0)$.
        \item Define $\hat{f}(i+1) := \max \{f(i) \cup \{last\}\}$ and put $last := \hat{f}(i+1)$.
    \end{enumerate}

    The morphism $\hat{f}$ is order-preserving. Since morphisms in $\Dalg^{op}$ are given by the preimages of morphisms in $\Dalg$, $\hat{f}(n) = m$ and therefore $\hat{f}$ is a cover.
\end{proof}

\begin{cor}
    The functor $\gamma: \Delta^{c} \to \Gamma$ is the opposite functor to the cardinality functor \hbox{$\card{-}:\Dalg\to\Fin\hookrightarrow\Fin_{\ast}$}.
\end{cor}

Consider the iterated wreath product $\Theta_n^{c} = \Delta^c \wr_{B} (\cdots \wr_{B}(\Delta^c \wr_{B} \Delta^c))$. The discussion above, together with Proposition \ref{prop-trees}, recovers the fact that the functor $\gamma_n: \Theta_n^{c} \to \Gamma$ is the opposite functor of the cardinality functor \hbox{$\card{-}: \Omega_n \to \Fin \hookrightarrow \Fin_{\ast}$\cite[Rem.~3.15.]{BERGER2007230}}. 

\subsection*{Wreath product of operator categories}
Operator categories, introduced by C. Barwick first in lecture notes \cite{barwick-notes} and later in an article \cite{barwick-paper}, can be viewed as an example of operadic categories. They are required to have a terminal object (so, in particular, they have only one connected component), and the fibers are given by pullbacks. To guarantee the axioms of operadic categories, a coherent choice of these pullbacks is required. Cardinalities are given by certain homsets. We refer the reader to the aforementioned references for details, with the remark that in \cite{barwick-notes} the author formulates the theory in the language of 1-categories and 1-operads, while in the later work \cite{barwick-paper} the machinery of $\infty$-categories and $\infty$-operads is employed. 

The wreath product of operator categories \cite[Ex.~1.6.]{barwick-paper} coincides with the wreath product of connected operadic categories. Our version of the wreath product is a natural generalization to the case when the categories in question have multiple connected components. The wreath product of operator categories is then used to introduce the \textit{external} Boardman-Vogt tensor product on the category of operads over various operator categories: given a $\Phi$-operad and $\opcat{P}$ and a $\Psi$-operad $\oper{Q}$, it produces a $(\Phi \wr \Psi)$-operad $\oper{P} _{\Phi}\otimes_{\Psi} \oper{Q}$. Similarly, in future work, we wish to introduce a kind of tensor product induced by the wreath product of operadic categories on the category $\Op$.

%% file: bibliography.bib
@article{opcat,
author = {Batanin, Michael and Markl, Martin},
year = {2014},
month = {04},
pages = {},
title = {Operadic categories and Duoidal {D}eligne's conjecture},
volume = {285},
journal = {Advances in Mathematics},
doi = {10.1016/j.aim.2015.07.008}
}

@article{Batanin-EH,
title = {The {E}ckmann–{H}ilton argument and higher operads},
journal = {Advances in Mathematics},
volume = {217},
number = {1},
pages = {334-385},
year = {2008},
issn = {0001-8708},
doi = {https://doi.org/10.1016/j.aim.2007.06.014},
url = {https://www.sciencedirect.com/science/article/pii/S0001870807002113},
author = {Michael Batanin}
}

@mastersthesis{masters,
    author = {Dunina, Daria},
    title = {Wreath product of operadic categories},
    school = {Charles University},
    year = {2024},
    url = {https://dspace.cuni.cz/handle/20.500.11956/193874}
}

@book{may,
    author = {May, Jon Peter},
    title = {The Geometry of Iterated Loop Spaces},
    publisher = {Springer Berlin, Heidelberg},
    year = {1972}
}

@article{grothendieck,
title = {The {G}rothendieck construction for model categories},
journal = {Advances in Mathematics},
volume = {281},
pages = {1306-1363},
year = {2015},
issn = {0001-8708},
doi = {https://doi.org/10.1016/j.aim.2015.03.031},
url = {https://www.sciencedirect.com/science/article/pii/S0001870815002030},
author = {Yonatan Harpaz and Matan Prasma}
}

@book{boardman-vogt,
    author = {Michael Boardman  and Rainer Vogt},
    title = {Homotopy Invariant Algebraic Structures on Topological Spaces},
    publisher = {Springer Berlin, Heidelberg},
    year = {1973}
}

@phdthesis{Brinkmeier2010,
author = {Michael Brinkmeier},
title = {On Operads},
year = {2010},
school = {Universit\"{a}t Osnabr\"{u}ck}
}

@ARTICLE{baratamoerdijk,
	author = {Barata, Miguel and Moerdijk, Ieke},
	title = {ON THE ADDITIVITY OF THE LITTLE CUBES OPERADS},
	year = {2025},
	journal = {Homology, Homotopy and Applications},
	volume = {27},
	number = {2},
	pages = {171 – 179},
	doi = {10.4310/HHA.2025.v27.n2.a8},
	type = {Article},
	publication_stage = {Final},
}

@article{BATANIN199839,
title = {Monoidal Globular Categories As a Natural Environment for the Theory of Weak n-Categories},
journal = {Advances in Mathematics},
volume = {136},
number = {1},
pages = {39-103},
year = {1998},
issn = {0001-8708},
doi = {https://doi.org/10.1006/aima.1998.1724},
url = {https://www.sciencedirect.com/science/article/pii/S0001870898917248},
author = {Michael Batanin}
}

@incollection {weiss,
    AUTHOR = {Weiss, Ittay},
     TITLE = {From operads to dendroidal sets},
 BOOKTITLE = {Mathematical foundations of quantum field theory and
              perturbative string theory},
    SERIES = {Proc. Sympos. Pure Math.},
    VOLUME = {83},
     PAGES = {31--70},
 PUBLISHER = {Amer. Math. Soc., Providence, RI},
      YEAR = {2011},
      ISBN = {978-0-8218-5195-1},
   MRCLASS = {55P48 (18D10 18D50 55U10 55U35)},
  MRNUMBER = {2742425},
MRREVIEWER = {Julia\ Bergner},
       DOI = {10.1090/pspum/083/2742425},
       URL = {https://doi.org/10.1090/pspum/083/2742425},
}

@article{Segal,
title = {Categories and cohomology theories},
journal = {Topology},
volume = {13},
number = {3},
pages = {293-312},
year = {1974},
issn = {0040-9383},
doi = {https://doi.org/10.1016/0040-9383(74)90022-6},
url = {https://www.sciencedirect.com/science/article/pii/0040938374900226},
author = {Graeme Segal}
}

@article{BERGER2007230,
title = {Iterated wreath product of the simplex category and iterated loop spaces},
journal = {Advances in Mathematics},
volume = {213},
number = {1},
pages = {230-270},
year = {2007},
issn = {0001-8708},
doi = {https://doi.org/10.1016/j.aim.2006.12.006},
url = {https://www.sciencedirect.com/science/article/pii/S0001870806003938},
author = {Clemens Berger}
}

@article{FIEDOROWICZ2015421,
title = {An Additivity Theorem for the interchange of En structures},
journal = {Advances in Mathematics},
volume = {273},
pages = {421-484},
year = {2015},
issn = {0001-8708},
doi = {https://doi.org/10.1016/j.aim.2014.10.020},
url = {https://www.sciencedirect.com/science/article/pii/S0001870814003661},
author = {Zbigniew Fiedorowicz and Rainer Vogt}
}

@article{brun,
    author = {Morten Brun and Zbigniew Fiedorowicz and Rainer Vogt},
    title = {On the multiplicative structure of topological Hochschild
homology},
    journal = {Algebraic \& Geometric Topology},
    year = 2007,
    volume = {7},
pages = {1633–1650}
}

@article{dunn,
title = {Tensor product of operads and iterated loop spaces},
journal = {Journal of Pure and Applied Algebra},
volume = {50},
number = {3},
pages = {237-258},
year = {1988},
issn = {0022-4049},
doi = {https://doi.org/10.1016/0022-4049(88)90103-X},
url = {https://www.sciencedirect.com/science/article/pii/002240498890103X},
author = {Gerald Dunn}
}

@unpublished{barwick-notes,
    author = {Clark Barwick},
    title = {Operator Categories, Multicategories, And Homotopy Coherent Algebra},
    note={lecture notes available at \url{https://webhomes.maths.ed.ac.uk/~cbarwick/papers/opcats.pdf}}
}

@article{barwick-paper,
author = {Barwick, Clark},
year = {2018},
month = {04},
pages = {1893-1959},
title = {From operator categories to higher operads},
volume = {22},
journal = {Geometry \& Topology},
doi = {10.2140/gt.2018.22.1893}
}

@article{trnka,
    author = {Dominik Trnka},
    title = {Integration of a categorical operad},
    journal = {Theory and Applications of Categories},
    volume = {45},
    pages = {130-150},
    number = {3},
    year = {2026}
}

@article{enriched-gr,
title = {The enriched Grothendieck construction},
journal = {Advances in Mathematics},
volume = {344},
pages = {234-261},
year = {2019},
issn = {0001-8708},
doi = {https://doi.org/10.1016/j.aim.2018.12.009},
url = {https://www.sciencedirect.com/science/article/pii/S0001870818305012},
author = {Jonathan Beardsley and Liang Ze Wong}
}
